\documentclass[oneside,french,american]{amsart}
\usepackage[T1]{fontenc}
\usepackage[latin9]{inputenc}
\usepackage{verbatim}
\usepackage{amsbsy}
\usepackage{amstext}
\usepackage{amsthm}
\usepackage{amssymb}
\usepackage{graphicx}
\PassOptionsToPackage{normalem}{ulem}
\usepackage{ulem}

\makeatletter

\providecommand{\tabularnewline}{\\}

\numberwithin{equation}{section}
\numberwithin{figure}{section}
\theoremstyle{plain}
\newtheorem{thm}{\protect\theoremname}[section]
\theoremstyle{plain}
\newtheorem{prop}[thm]{\protect\propositionname}
\theoremstyle{plain}
\newtheorem{fact}[thm]{\protect\factname}
\theoremstyle{plain}
\newtheorem{lem}[thm]{\protect\lemmaname}
\theoremstyle{plain}
\newtheorem{cor}[thm]{\protect\corollaryname}
\theoremstyle{remark}
\newtheorem{rem}[thm]{\protect\remarkname}

\makeatletter

\usepackage{bbm}
\usepackage{graphicx}
\usepackage{tikz-qtree}

\allowdisplaybreaks

\newtheorem{hypothesis}{Assumption}

\newtheorem{subhyp}{Assumption}

\numberwithin{equation}{section}

\newcommand{\1}{\mathbbm{1}}

\newcommand{\N}{\mathbb{N}}
\newcommand{\Z}{\mathbb{Z}}

\newcommand{\E}{\mathbb{E}}
\newcommand{\R}{\mathbb{R}}
\newcommand{\p}{\mathbb{P}}

\renewcommand{\epsilon}{ \varepsilon}

\newcommand{\Bsym}{\mathcal{B}_{\text{sym}}^0}

\newcommand{\Ue}{\mathcal{U}_\epsilon}

\DeclareMathOperator*{\Id}{Id}

\makeatother

\usepackage{a4wide}

\makeatother

\usepackage{babel}
\makeatletter
\addto\extrasfrench{%
   \providecommand{\fg}{\ifdim\lastskip>\z@\unskip\fi~\frqq}%
}

\makeatother
\addto\captionsamerican{\renewcommand{\corollaryname}{Corollary}}
\addto\captionsamerican{\renewcommand{\factname}{Fact}}
\addto\captionsamerican{\renewcommand{\lemmaname}{Lemma}}
\addto\captionsamerican{\renewcommand{\propositionname}{Proposition}}
\addto\captionsamerican{\renewcommand{\remarkname}{Remark}}
\addto\captionsamerican{\renewcommand{\theoremname}{Theorem}}
\addto\captionsfrench{\renewcommand{\corollaryname}{Corollaire}}
\addto\captionsfrench{\renewcommand{\factname}{Fait}}
\addto\captionsfrench{\renewcommand{\lemmaname}{Lemme}}
\addto\captionsfrench{\renewcommand{\propositionname}{Proposition}}
\addto\captionsfrench{\renewcommand{\remarkname}{Remarque}}
\addto\captionsfrench{\renewcommand{\theoremname}{Théorème}}
\providecommand{\corollaryname}{Corollary}
\providecommand{\factname}{Fact}
\providecommand{\lemmaname}{Lemma}
\providecommand{\propositionname}{Proposition}
\providecommand{\remarkname}{Remark}
\providecommand{\theoremname}{Theorem}

\begin{document}

\title{A central-limit Theorem for conservative fragmentation chains}

\author{Camille Noûs, Sylvain Rubenthaler,}

\email{camille.nous@cogitamus.fr, rubentha@unice.fr}

\date{2022}
\begin{abstract}
We are interested in a fragmentation process. We observe fragments
frozen when their sizes are less than $\epsilon$ ($\epsilon>0$).
It is known (\cite{bertoin-martinez-2005}) that the empirical measure
of these fragments converges in law, under some renormalization. In
\cite{hoffmann-krell-2011}, the authors show a bound for the rate
of convergence. Here, we show a central-limit theorem, under some
assumptions. This gives us an exact rate of convergence. 
\end{abstract}


\maketitle

\section{Introduction}

\subsection{Scientific and economic context}

One of the main goals in the mining industry is to extract blocks
of metallic ore and then separate the metal from the valueless material.
To do so, rock is fragmented into smaller and smaller rocks. This
is carried out in a series of steps, the first one being blasting,
after which the material goes through a sequence of crushers. At each
step, the particles are screened, and if they are smaller than the
diameter of the mesh of a classifying grid, they go to the next crusher.
The process stops when the material has a sufficiently small size
(more precisely, small enough to enable physicochemical processing).

This fragmentation process is energetically costly (each crusher consumes
a certain quantity of energy to crush the material it is fed). One
of the problems that faces the mining industry is that of minimizing
the energy used. The optimization parameters are the number of crushers
and the technical specifications of these crushers.

In \cite{bertoin-martinez-2005}, the authors propose a mathematical
model of what happens in a crusher. In this model, the rock pieces/fragments
are fragmented independently of each other, in a random and auto-similar
manner. This is consistent with what is observed in the industry,
and this is supported by the following publications: \cite{bird-perrier-2002,devoto-martinez-1998,weiss-1986,turcotte-1986}.
Each fragment has a size $s$ (in $\R^{+}$) and is then fragmented
into smaller fragments of sizes $s_{1}$, $s_{2}$, \ldots{} such that
the sequence $(s_{1}/s,s_{2}/s,\dots)$ has a law $\nu$ which does
not depend on $s$ (which is why the fragmentation is said to be auto-similar).
This law $\nu$ is called the \emph{dislocation measure} (each crusher
has its own dislocation measure). The dynamic of the fragmentation
process is thus modeled in a stochastic way.

In each crusher, the rock pieces are fragmented repetitively until
they are small enough to slide through a mesh whose holes have a fixed
diameter. So the fragmentation process stops for each fragment when
its size is smaller than the diameter of the mesh, which we denote
by $\epsilon$ ($\epsilon>0$). We are interested in the \emph{statistical
distribution }of the fragments coming out of a crusher. If we renormalize
the sizes of these fragments by dividing them by $\epsilon$, we obtain
a measure $\gamma_{-\log(\epsilon)}$, which we call the \emph{empirical
measure }(the reason for the index $-\log(\epsilon)$ instead of $\epsilon$
will be made clear later). In \cite{bertoin-martinez-2005}, the authors
show that the energy consumed by the crusher to reduce the rock pieces
to fragments whose diameters are smaller than $\epsilon$ can be computed
as an integral of a bounded function against the measure $\gamma_{-\log(\epsilon)}$
(they cite \cite{bond-1952,charles-1957,walker-lewis-mcadames-gilliland-1985}
on this particular subject). For each crusher, the empirical measure
$\gamma_{-\log(\epsilon)}$ is one of the two only observable variables
(the other one being the size of the pieces pushed into the grinder).
The specifications of a crusher are summarized in $\epsilon$ and
$\nu$.

\subsection{State of the art}

In \cite{bertoin-martinez-2005}, the authors show that the energy
consumed by a crusher to reduce rock pieces of a fixed size into fragments
whose diameter are smaller than $\epsilon$ behaves asymptotically
like a power of $\epsilon$ when $\epsilon$ goes to zero. More precisely,
this energy multiplied by a power of $\epsilon$ converges towards
a constant of the form $\kappa=\nu(\varphi)$ (the integral of $\nu$,
the dislocation measure, against a bounded function $\varphi$). In
\cite{bertoin-martinez-2005}, the authors also show a law of large
numbers for the empirical measure $\gamma_{-\log(\epsilon)}$. More
precisely, if $f$ is bounded continuous, $\gamma_{-\log(\epsilon)}(f)$
converges in law, when $\epsilon$ goes to zero, towards an integral
of $f$ against a measure related to $\nu$ (this result also appears
in \cite{hoffmann-krell-2011}, p. 399). We set $\gamma_{\infty}(f)$
to be this limit (check Equations (\ref{eq:def-gamma-infty}), (\ref{eq:def-eta}),
(\ref{eq:def-pi}) to get an exact formula). The empirical measure
$\gamma_{-\log(\epsilon)}$ thus contains information relative to
$\nu$ and one could extract from it an estimation of $\kappa$ or
of an integral of any function against $\nu$. 

It is worth noting that by studying what happens in various crushers,
we could study a family $(\nu_{i}(f_{j}))_{i\in I,j\in J}$ (with
an index $i$ for the number of the crusher and the index $j$ for
the $j$-th test function in a well-chosen basis). Using statistical
learning methods, one could from there make a prediction for $\nu(f_{j}$)
for a new crusher for which we know only the mechanical specifications
(shape, power, frequencies of the rotating parts \ldots ). It would
evidently be interesting to know $\nu$ before even building the crusher.

In the same spirit, \cite{fontbona-krell-martinez-2010} studies the
energy efficiency of two crushers used after one another. When the
final size of the fragments tends to zero, this paper tells us wether
it is more efficient enrgywise to use one crusher or two crushers
in a row (another asymptotic is also considered in the paper). 

In \cite{harris-knobloch-kyprianou-2010}, the authors prove a convergence
result for the empirical measure similar to the one in \cite{bertoin-martinez-2005},
the convergence in law being replaced by an almost sure convergence.
In \cite{hoffmann-krell-2011}, the authors give a bound on the rate
of this convergence, in a $L^{2}$ sense, under the assumption that
the fragmentation is conservative. This assumption means there is
no loss of mass due to the formation of dust during the fragmentation
process.\foreignlanguage{french}{}
\begin{figure}[h]
\selectlanguage{french}%
\[
\begin{array}{ccc}
\gamma_{-\log(\epsilon)} & \underset{\underset{\epsilon\rightarrow0}{\longrightarrow}}{\text{(bound on rate)}} & \gamma_{\infty}\\
 &  & \updownarrow\text{relation }\\
\text{energy}\times(\text{power of \ensuremath{\epsilon})} & \underset{\epsilon\rightarrow0}{\sim} & \kappa=\nu(\varphi)
\end{array}
\]

\caption{\foreignlanguage{american}{State of the art.}}
\selectlanguage{american}%
\end{figure}

So we have convergence results (\cite{bertoin-martinez-2005,harris-knobloch-kyprianou-2010})
of an empirical quantity towards constants of interest (a different
constant for each test function $f$). Using some transformations,
these constants could be used to estimate the constant $\kappa$.
Thus it is natural to ask what is the exact rate of convergence in
this estimation, if only to be able to build confidence intervals.
In \cite{hoffmann-krell-2011}, we only have a bound on the rate.

When a sequence of empirical measures converges to some measure, it
is natural to study the fluctuations, which often turn out to be Gaussian.
For such results in the case of empirical measures related to the
mollified Boltzmann equation, one can cite \cite{meleard-1998,uchiyama-1988,dawson-zheng-1991}.
When interested in the limit of a $n$-tuple as in Equation (\ref{eq:n-uplet})
below, we say we are looking at the convergence of a $U$-statistics.
Textbooks deal with the case where the points defining the empirical
measure are independent or with a known correlation (see \cite{delapena-gine-1999,dynkin-mandelbaum-1983,lee-1990}).
The problem is more complex when the points defining the empirical
measure are in interaction with each other like it is the case here.

\subsection{Goal of the paper}

As explained above, we want to obtain the rate of convergence in the
convergence of $\gamma_{-\log(\epsilon)}$ when $\epsilon$ goes to
zero. We want to produce a central-limit theorem of the kind: for
a bounded continuous $f$, $\epsilon^{\beta}(\gamma_{-\log(\epsilon)}(f)-\gamma_{\infty}(f))$
converges towards a non-trivial measure when $\epsilon$ goes to zero
(the limiting measure will in fact be Gaussian), for some exponent
$\beta$. The technics used will allow us to prove the convergence
towards a multivariate Gaussian of a vector of the kind 
\begin{equation}
\epsilon^{\beta}(\gamma_{-\log(\epsilon)}(f_{1})-\gamma_{\infty}(f_{1}),\dots,\gamma_{-\log(\epsilon)}(f_{n})-\gamma_{\infty}(f_{n}))\label{eq:n-uplet}
\end{equation}
for functions $f_{1}$, \ldots , $f_{n}$.

More precisely, if by $Z_{1}$, $Z_{2}$, \ldots , $Z_{N}$ we denote
the fragments sizes that go out from a crusher (with mesh diameter
equal to $\epsilon$). We would like to show that for a bounded continuous
$f$, 
\[
\gamma_{-\log(\epsilon)}(f):=\sum_{i=1}^{N}Z_{i}f(Z_{i}/\epsilon)\longrightarrow\gamma_{\infty}(f)\text{, almost surely, when }\epsilon\rightarrow0\,,
\]
and that for all $n$, and $f_{1}$, \ldots ,$f_{n}$ bounded continuous
function such that $\gamma_{\infty}(f_{i})=0$,
\[
\epsilon^{\beta}(\gamma_{-\log(\epsilon)}(f_{1}),\dots,\gamma_{-\log(\epsilon)}(f_{n}))
\]
converges in law towards a multivariate Gaussian when $\epsilon$
goes to zero. 

The exact results are stated in Proposition \ref{prop:convergence-ps}
and Theorem \ref{thm:central-limit}.

\subsection{Outline of the paper}

We will state our assumptions along the way (Assumptions \ref{hyp:A},
\ref{hyp:conservative}, \ref{hyp:delta-step}, \ref{hyp:queue-pi}).
Assumption \ref{hyp:queue-pi} can be found at the beginning of Section
\ref{sec:Rate-of-convergence}. We define our model in Section \ref{sec:Statistical-model}.
The main idea is that we want to follow tags during the fragmentation
process. Let us imagine the fragmentation is the process of breaking
a stick (modeled by $[0,1]$) into smaller sticks. We suppose that
the original stick has painted dots and that during the fragmentation
process, we take note of the sizes of the sticks supporting the painted
dots. When the sizes of these sticks get smaller than $\epsilon$
($\epsilon>0$), the fragmentation is stopped for them and we call
them the painted sticks. In Section \ref{sec:Rate-of-convergence},
we make use of classical results on renewal processes and of \cite{sgibnev-2002}
to show that the size of one painted stick has an asymptotic behavior
when $\epsilon$ goes to zero and that we have a bound on the rate
with which it reaches this behavior. Section \ref{sec:Limits-of-symmetric}
is the most technical. There we study the asymptotics of symmetric
functionals of the sizes of the painted sticks (always when $\epsilon$
goes to zero). In Section \ref{sec:Results}, we precisely define
the measure we are interested in ($\gamma_{T}$ with $T=-\log(\epsilon)$).
Using the results of Section \ref{sec:Limits-of-symmetric}, it is
then easy to show a law of large numbers for $\gamma_{T}$ (Proposition
\ref{prop:convergence-ps}) and a central-limit Theorem (Theorem \ref{thm:central-limit}).
Proposition \ref{prop:convergence-ps} and Theorem \ref{thm:central-limit}
are our two main results. The proof of Theorem \ref{thm:central-limit}
is based on a simple computation involving characteristic functions
(the same technique was already used in \cite{del-moral-patras-rubenthaler-2009,del-moral-patras-rubenthaler-2011a,del-moral-patras-rubenthaler-2011b,rubenthaler-2016}).

\subsection{Notations\label{subsec:Notations}}

For $x$ in $\R$, we set $\lceil x\rceil=\inf\{n\in\Z\,:\,n\geq x\}$,
$\lfloor x\rfloor=\sup\{n\in\Z\,:\,n\leq x\}$. The symbol $\sqcup$
means ``disjoint union''. For $n$ in $\N^{*}$, we set $[n]=\{1,2,\dots,n\}$.
For $f$ an application from a set $E$ to a set $F$, we write $f:E\hookrightarrow F$
if $f$ is injective and, for $k$ in $\N^{*}$, if $F=E$, we set
\[
f^{\circ k}=\underset{k\mbox{ times }}{\underbrace{f\circ f\circ\dots\circ f}}\,.
\]
For any set $E$, we set $\mathcal{P}(E)$ to be the set of subsets
of $E$. 

\section{Statistical model\label{sec:Statistical-model}}

\subsection{Fragmentation chains\label{subsec:Fragmentation-chains}}

Let $\epsilon>0$. Like in \cite{hoffmann-krell-2011}, we start with
the space 
\[
\mathcal{S}^{\downarrow}=\left\{ \mathbf{s}=(s_{1},s_{2},\dots),\,s_{1}\geq s_{2}\geq\dots\geq0,\,\sum_{i=1}^{+\infty}s_{i}\leq1\right\} \,.
\]
A fragmentation chain is a process in $\mathcal{S}^{\downarrow}$
characterized by 
\begin{itemize}
\item a dislocation measure $\nu$ which is a finite measure on $\mathcal{S}^{\downarrow}$,
\item a description of the law of the times between fragmentations.
\end{itemize}
A fragmentation chain with dislocation measure $\nu$ is a Markov
process $X=(X(t),t\geq0)$ with values in $\mathcal{S}^{\downarrow}$.
Its evolution can be described as follows: a fragment with size $x$
lives for some time (which may or may not be random) then splits and
gives rise to a family of smaller fragments distributed as $x\xi$,
where $\xi$ is distributed according to $\nu(.)/\nu(\mathcal{S}^{\downarrow})$.
We suppose the life-time of a fragment of size $x$ is an exponential
time of parameter $x^{\alpha}\nu(\mathcal{S}^{\downarrow})$, for
some $\alpha$. We could here make different assumptions on the life-time
of fragments, but this would not change our results. Indeed, as we
are interested in the sizes of the fragments frozen as soon as they
are smaller than $\epsilon$, the time they need to become this small
is not important. 

We denote by $\p_{m}$ the law of $X$ started from the initial configuration
$(m,0,0,\dots)$ with $m$ in $(0,1]$. The law of $X$ is entirely
determined by $\alpha$ and $\nu(.)$ (Theorem 3 of \cite{bertoin-2002}). 

We make the same assumption as in \cite{hoffmann-krell-2011} and
we will call it Assumption \ref{hyp:A}.

\begin{hypothesis}\label{hyp:A}

We have $\nu(\mathcal{S}^{\downarrow})=1$ and $\nu(s_{1}\in]0;1[)=1$.

\end{hypothesis}

Let
\[
\mathcal{U}:=\{0\}\cup\bigcup_{n=1}^{+\infty}(\N^{*})^{n}
\]
denote the infinite genealogical tree. For $u$ in $\mathcal{U}$,
we use the conventional notation $u=()$ if $u=\{0\}$ and $u=(u_{1},\dots,u_{n})$
if $u\in(\N^{*})^{n}$ with $n\in\N^{*}.$ This way, any $u$ in $\mathcal{U}$
can be denoted by $u=(u_{1},\dots,u_{n})$, for some $u_{1},\dots,u_{n}$
and with $n$ in $\N$. Now, for $u=(u_{1},\dots,u_{n})\in\mathcal{U}$
and $i\in\N^{*}$, we say that $u$ is in the $n$-th generation and
we write $|u|=n$, and we write $ui=(u_{1},\dots,u_{n},i)$, $u(k)=(u_{1},\dots,u_{k})$
for all $k\in[n]$. %
{}  For any $u=(u_{1},\dots,u_{n})$ and $v=ui$ ($i\in\N^{*}$), we
say that $u$ is the mother of $v$. For any $u$ in $\mathcal{U}\backslash\{0\}$
($\mathcal{U}$ deprived of its root), $u$ has exactly one mother
and we denote it by $\boldsymbol{m}(u)$. The set $\mathcal{U}$ is
ordered alphanumerically~:
\begin{itemize}
\item If $u$ and $v$ are in $\mathcal{U}$ and $|u|<|v|$ then $u<v$.
\item If $u$ and $v$ are in $\mathcal{U}$ and $|u|=|v|=n$ and $u=(u_{1},\dots,u_{n})$,
$v=(v_{1},\dots,v_{n})$ with $u_{1}=v_{1}$, \dots{} , $u_{k}=v_{k}$,
$u_{k+1}<v_{k+1}$ then $u<v$.
\end{itemize}
Suppose we have a process $X$ which has the law $\p_{m}$. For all
$\omega$, we can index the fragments that are formed by the process
$X$ with elements of $\mathcal{U}$ in a recursive way.
\begin{itemize}
\item We start with a fragment of size $m$ indexed by $u=()$. 
\item If a fragment $x$, with a birth-time $t_{1}$ and a split-time $t_{2}$,
is indexed by $u$ in $\mathcal{U}$. At time $t_{2}$, this fragment
splits into smaller fragments of sizes $(xs_{1},xs_{2},\dots)$ with
$(s_{1},s_{2},\dots)$ of law $\nu(.)/\nu(\mathcal{S}^{\downarrow})$.
We index the fragment of size $xs_{1}$ by $u1$, we index the fragment
of size $xs_{2}$ by $u2$, and so on. 
\end{itemize}
A mark is an application from $\mathcal{U}$ to some other set. We
associate a mark $\xi_{\dots}$ on the tree $\mathcal{U}$ to each
path of the process $X$. The mark at node $u$ is $\xi_{u}$, where
$\xi_{u}$ is the size of the fragment indexed by $u$.%
{} The distribution of this random mark can be described recursively
as follows.
\begin{prop}
(Consequence \label{prop:(reformulation-of-Proposition} of Proposition
1.3, p. 25, \cite{bertoin-2006}) There exists a family of i.i.d.
variables indexed by the nodes of the genealogical tree, $((\widetilde{\xi}_{ui})_{i\in\N^{*}},u\in\mathcal{U})$,
where each $(\widetilde{\xi}_{ui})_{i\in\N^{*}}$ is distributed according
to the law $\nu(.)/\nu(\mathcal{S}^{\downarrow})$, and such that
the following holds:\\
Given the marks $(\xi_{v},|v|\leq n)$ of the first $n$ generations,
the marks at generation $n+1$ are given by 
\[
\xi_{ui}=\widetilde{\xi}_{ui}\xi_{u}\,,
\]
where $u=(u_{1,}\dots,u_{n})$ and $ui=(u_{1},\dots,u_{n},i)$ is
the $i-th$ child of $u$.
\end{prop}

\subsection{Tagged fragments\label{subsec:Tagged-fragments}}

From now on, we suppose that we start with a block of size $m=1$.
We assume that the total mass of the fragments remains constant through
time, as follows.

\begin{hypothesis}\textbf{\label{hyp:conservative}(Conservative
property).}

We have $\nu(\sum_{i=1}^{+\infty}s_{i}=1)=1$.

\end{hypothesis}

This assumption was already present in \cite{hoffmann-krell-2011}.
We observe that the Malthusian exponent of \cite{bertoin-2006} (p.
27) is equal to $1$ under our assumptions. Without this assumption,
the link between the empirical measure $\gamma_{-\log(\epsilon)}$
and the tagged fragments (Equation (\ref{eq:lien-odot-B})) vanishes
and our proofs of Proposition \ref{prop:convergence-ps} and Theorem
\ref{thm:central-limit} fail.

We can now define tagged fragments.%
{} We use the representation of fragmentation chains as random infinite
marked tree to define a fragmentation chain with $q$ tags. Suppose
we have a fragmentation process $X$ of law $\p_{1}$. We take $(Y_{1},Y_{2},\dots,Y_{q})$
to be $q$ i.i.d. variables of law $\mathcal{U}([0,1])$. We set,
for all $u$ in $\mathcal{U}$,
\[
(\xi_{u},A_{u},I_{u})
\]
with $\xi_{u}$ defined as above. The random variables $A_{u}$ take
values in the subsets of $[q]$. The random variables $I_{u}$ are
intervals. These variables are defined as follows.
\begin{itemize}
\item We set $A_{\{0\}}=[q]$, $I_{\{0\}}=(0,1]$ ($I_{\{0\}}$ is of length
$\xi_{\{0\}}=1$)
\item For all $n\in\N$. Suppose we are given the marks of the first $n$
generations. Suppose that, for $u$ in the $n$-th generation, $I_{u}=(a_{u},a_{u}+\xi_{u}]$
for some $a_{u}\in\R$ (it is of length $\xi_{u}$). Then the marks
at generation $n+1$ are given by Proposition \ref{prop:(reformulation-of-Proposition}
(concerning $\xi_{.}$) and, for all $u$ such that $|u|=n$ and for
all $i$ in $\N^{*}$ 
\[
I_{ui}=(a_{u}+\xi_{u}(\widetilde{\xi}_{u1}+\dots+\widetilde{\xi}_{u(i-1)}),a_{u}+\xi_{u}(\widetilde{\xi}_{u1}+\dots+\widetilde{\xi}_{ui})]\,,
\]
\[
k\in A_{ui}\text{ if and only if }Y_{k}\in I_{ui}\,,
\]
($I_{ui}$ is then of length $\xi_{ui}$). We observe that for all
$j\in[q]$, $u\in\mathcal{U}$, $i\in\N^{*}$, 
\begin{equation}
\p(j\in A_{ui}|j\in A_{u},\widetilde{\xi}_{ui})=\widetilde{\xi}_{ui}\,.\label{eq:proba-rester-dans-fragment}
\end{equation}
\end{itemize}
In this definition, we imagine having $q$ dots on the interval $[0,1]$
and we impose that the dot $j$ has the position $Y_{j}$ (for all
$j$ in $[q]$). During the fragmentation process, if we know that
the dot $j$ is in the interval $I_{u}$ of length $\xi_{u}$, then
the probability that this dot is on $I_{ui}$ (for some $i$ in $\N^{*}$,
$I_{ui}$ of length $\xi_{ui}$) is equal to $\xi_{ui}/\xi_{u}=\widetilde{\xi}_{ui}$. 

In the case $q=1$, the branch $\{u\in\mathcal{U}\,:\,A_{u}\neq\emptyset\}$
has the same law as the randomly tagged branch of Section 1.2.3 of
\cite{bertoin-2006}. The presentation is simpler in our case because
the Malthusian exponent is $1$ under Assumption \ref{hyp:conservative}.

\subsection{Observation scheme\label{subsec:Observation-scheme}}

We freeze the process when the fragments become smaller than a given
threshold $\epsilon>0$. That is, we have the following data
\[
(\xi_{u})_{u\in\mathcal{U}_{\epsilon}}\,,
\]
where 
\[
\mathcal{U}_{\epsilon}=\{u\in\mathcal{U},\,\xi_{\boldsymbol{m}(u)}\geq\epsilon,\,\xi_{u}<\epsilon\}\,.
\]

We now look at $q$ tagged fragments ($q\in\N^{*}$). For each $i$
in $[q]$, we call $L_{0}^{(i)}=1$, $L_{1}^{(i)}$, $L_{2}^{(i)}$\dots{}
the successive sizes of the fragment having the tag $i$. More precisely,
for each $n\in\N^{*}$, there is almost surely exactly one $u\in\mathcal{U}$
such that $|u|=n$ and $i\in A_{u}$; and so, $L_{n}^{(i)}=\xi_{u}$.
For each $i$, the process $S_{0}^{(i)}=-\log(L_{0}^{(i)})=0\leq S_{1}^{(i)}=-\log(L_{1}^{(i)})\leq\dots$
is a renewal process without delay, with waiting-time following a
law $\pi$ (see \cite{asmussen-2003}, Chapter V for an introduction
to renewal processes). The waiting times are (for $i$ in $[q]$):
$S_{0}^{(i)}$, $S_{1}^{(i)}-S_{0}^{(i)}$, $S_{2}^{(i)}-S_{1}^{(i)}$,
\dots{} The renewal times are (for $i$ in $[q]$): $S_{0}^{(i)}$,
$S_{1}^{(i)}$, $S_{2}^{(i)}$, \dots{} The law $\pi$ is defined
by the following.
\begin{equation}
\text{For all bounded measurable }f:[0,1]\rightarrow[0,+\infty)\,,\,\int_{\mathcal{S}^{\downarrow}}\sum_{i=1}^{+\infty}s_{i}f(s_{i})\nu(d\boldsymbol{s})=\int_{0}^{+\infty}f(e^{-x})\pi(dx)\,,\label{eq:def-pi}
\end{equation}
(see Proposition 1.6, p. 34 of \cite{bertoin-2006}, or Equations
(3), (4), p. 398 of \cite{hoffmann-krell-2011}). Under Assumption
\ref{hyp:A} and Assumption \ref{hyp:conservative}, Proposition 1.6
of \cite{bertoin-2006} is true, even without the Malthusian Hypothesis
of \cite{bertoin-2006}.

We make the following assumption on $\pi$. 

\stepcounter{hypothesis}

\begin{subhyp}\label{hyp:delta-step}

There exist $a$ and $b$ ($0<a<b<+\infty$) such that the support
of $\pi$ is $[a,b]$. We set $\delta=e^{-b}$. 

\end{subhyp}

We added a comment about the above Assumption in Remark \ref{rem:About-Assumption}.
We believe that we could replace the above Assumption by the following.

\begin{subhyp}\label{hyp:no-step}

The support of $\pi$ is $(0,+\infty)$.

\end{subhyp}

However, this would lead to difficult computations.

We set 
\begin{equation}
T=-\log(\epsilon)\,.\label{eq:relation-T-epsilon}
\end{equation}
 We set, for all $i\in[q]$, $t\geq0$,
\[
N_{t}^{(i)}=\inf\{j\,:\,S_{j}^{(i)}>t\}\,,
\]
 \textbf{
\begin{equation}
B_{t}^{(i)}=S_{N_{t}^{(i)}}^{(i)}-t\,,\label{eq:def-B}
\end{equation}
\[
C_{t}^{(i)}=t-S_{N_{t}^{(i)}-1}^{(i)}\,,
\]
}
\begin{figure}[h]
\centering{}\includegraphics[scale=0.7]{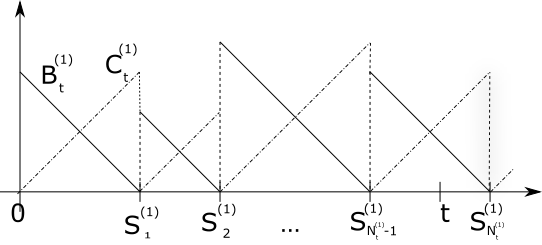}\textbf{\caption{Process\textbf{ $B^{(1)}$ and $C^{(1)}$.\label{fig:Process-.}}}
}
\end{figure}
(see Figure \ref{fig:Process-.} for an illustration). The processes
$B^{(i)}$, $C^{(i)}$, $N^{(i)}$ are time-homogeneous Markov processes
(Proposition 1.5 p. 141 of \cite{asmussen-2003}). All of them are
càdlàg (i.e. right-continuous with a left-hand side limit). We call
$B^{(i)}$ the residual lifetime of the fragment tagged by $i$. We
call $C^{(i)}$ the age of the fragment tagged by $i$. We call $N^{(i)}$
the number of renewal before $t$. In the following, we will treat
$t$ as a time parameter. This has nothing to do with the time in
which the fragmentation process $X$ evolves. 

We observe that, for all $t$, $(B_{t}^{(1)},\dots,B_{t}^{(q)})$
is exchangeable (meaning that for all $\sigma$ in the symmetric group
of order $q$, $(B_{t}^{(\sigma(1))},\dots,B_{t}^{(\sigma(q))})$
has the same law as $(B_{t}^{(1)},\dots,B_{t}^{(q)})$).

When we look at the fragments of sizes $(\xi_{u},\,u\in\mathcal{U}_{\epsilon}\,:\,A_{u}\neq\emptyset)$,
we have almost the same information as when we look at $(B_{T}^{(1)},B_{T}^{(2)},\dots,B_{T}^{(q)})$.
We say \emph{almost} because knowing $(B_{T}^{(1)},B_{T}^{(2)},\dots,B_{T}^{(q)})$
does not give exactly the number of $u$ in $\mathcal{U}_{\epsilon}$
such that $A_{u}$ is not empty.

In the remaining of Section \ref{sec:Statistical-model}, we define
processes that will be useful when we will describe the asymptotics
of our model (in Section \ref{sec:Limits-of-symmetric}). 

\subsection{Stationary renewal processes ($\overline{B}^{(1)}$, $\overline{B}^{(1),v}$).
\label{subsec:Stationary-age-process}}

We define $\widetilde{X}$ to be an independent copy of $X$. We suppose
it has $q$ tagged fragments. Therefore it has a mark $(\widetilde{\xi},\widetilde{A})$
and a renewal processes $(\widetilde{S}_{k}^{(i)})_{k\geq0}$ (for
all $i$ in $[q]$) defined in the same way as for $X$. We let $(\widetilde{B}^{(1)},\widetilde{B}^{(2)})$
be the residual lifetimes of the fragments tagged by $1$ and $2$.

Let 
\[
\mu=\int_{0}^{+\infty}x\pi(dx)
\]
 and let $\pi_{1}$ be the distribution with density $x\mapsto x/\mu$
with respect to $\pi$. We set $\overline{C}$ to be a random variable
of law $\pi_{1}$. We set $U$ to be independent of $\overline{C}$
and uniform on $(0,1)$. We set $\widetilde{S}_{-1}=\overline{C}(1-U)$.
The process $\overline{S}_{0}=\widetilde{S}_{-1}$, $\overline{S}_{1}=\widetilde{S}_{-1}+\widetilde{S}_{0}^{(1)}$
, $\overline{S}_{2}=\widetilde{S}_{-1}+\widetilde{S}_{1}^{(1)}$,
$\overline{S}_{2}=\widetilde{S}_{-1}+\widetilde{S}_{2}^{(1)}$, \ldots{}
is a renewal process with delay $\pi_{1}$ (with waiting times $\overline{S}_{0}$,
$\overline{S}_{1}-\overline{S}_{0}$, \dots{} all smaller than $b$
by Assumption \ref{hyp:delta-step}). The renewal times are $\overline{S}_{0}$,
$\overline{S}_{1}$, $\overline{S}_{2}$, \dots{} We set $(\overline{B}_{t}^{(1)})_{t\geq0}$
to be the residual lifetime process of this renewal process:
\begin{equation}
\overline{B}_{t}^{(1)}=\begin{cases}
\overline{C}(1-U)-t & \mbox{ if }t<\overline{S}_{0}\,,\\
\inf_{n\geq0}\{\overline{S}_{n}\,:\,\overline{S}_{n}>t\}-t & \mbox{ if }t\geq\overline{S}_{0}\,,
\end{cases}\label{eq:def-B-stationnaire}
\end{equation}
and we define $(\overline{C}_{t}^{(1)})_{t\geq0}$:
\begin{equation}
\overline{C}_{t}^{(1)}=\begin{cases}
\overline{C}U+t & \mbox{ if }t<\overline{S}_{0}\,,\\
t-\sup_{n\geq0}\{\overline{S}_{n}\,:\,\overline{S}_{n}\leq t\} & \mbox{ if }t\geq\overline{S}_{0}\,,
\end{cases}\label{eq:def-C-stationnaire}
\end{equation}
(we call it the age process of our renewal process) and we set
\[
\overline{N}_{t}^{(1)}=\inf\{j\,:\,\overline{S}_{j}>t\}\,.
\]

\begin{fact}
\label{fact:Theorem-3.3-p.151}Theorem 3.3 p.151 of \cite{asmussen-2003}
tells us that $(\overline{B}_{t}^{(1)},\overline{C}_{t}^{(1)})_{t\geq0}$
has the same transition as $(B_{t}^{(1)},C_{t}^{(1)})_{t\geq0}$ defined
above and that $(\overline{B}_{t}^{(1)},\overline{C}_{t}^{(1)})_{t\geq0}$
is stationary. In particular, this means that the law of $\overline{B}_{t}^{(1)}$
does not depend on $t$.
\end{fact}

\begin{figure}[h]
\begin{centering}
\includegraphics[scale=0.75]{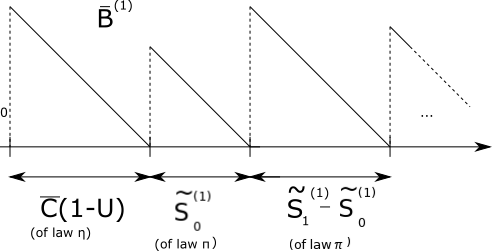}\caption{Renewal process with delay\label{fig:Renewal-process-with}}
\par\end{centering}
\end{figure}

We see in Figure \ref{fig:Renewal-process-with} a graphic representation
of $\overline{B}_{.}^{(1)}$. 
This might be counter-intuitive to start with $\overline{B}_{0}^{(1)}$having
a law which is not $\pi$ in order to get a stationary process, but
Corollary 3.6 p. 153 of \cite{asmussen-2003} is clear on this point:
a delayed renewal process (with waiting-time of law $\pi$) is stationary
if and only if the distribution of the initial delay is $\eta$ (defined
below). 

We define a measure $\eta$ on $\R^{+}$ by its action on bounded
measurable functions:
\begin{equation}
\text{For all bounded measurable }f\,:\,\R^{+}\rightarrow\R\,,\,\eta(f)=\frac{1}{\mu}\int_{\R^{+}}\E(f(Y-s)\1_{\{Y-s>0\}})ds\,,\,(Y\sim\pi)\,.\label{eq:def-eta}
\end{equation}

\begin{lem}
\label{lem:The-measure-}The measure $\eta$ is the law of $\overline{B}_{t}^{(1)}$
(for any $t\geq0$). It is also the law of $(\overline{C}_{t}^{(1)})$
(for any $t$). 
\end{lem}

\begin{proof}
We write the proof for $\overline{B}_{t}^{(1)}$only. Let $\xi\geq0$.
We set $f(y)=\1_{y\geq\xi}$, for all $y$ in $\R$. We have (with
$Y$ of law $\pi$)
\begin{eqnarray*}
\frac{1}{\mu}\int_{\R^{+}}\E(f(Y-s)\1_{Y-s>0})ds & = & \frac{1}{\mu}\int_{\R^{+}}\left(\int_{0}^{y}\1_{y-s\geq\xi}ds\right)\pi(dy)\\
 & = & \frac{1}{\mu}\int_{\R^{+}}(y-\xi)_{+}\pi(dy)\\
 & = & \int_{\xi}^{+\infty}\left(1-\frac{\xi}{y}\right)\frac{y}{\mu}\pi(dy)\\
 & = & \p(\overline{C}(1-U)\geq\xi)\,.
\end{eqnarray*}
\end{proof}
We set $\eta_{2}$ to be the law of $(\overline{C}_{0}^{(1)},\overline{B}_{0}^{(1)})=(\overline{C}U,\overline{C}(1-U)).$
The support of $\eta_{2}$ is $\mathcal{C}:=\{(u,v)\in[0,b]^{2}\,:\,a\leq u+v\leq b\}$.

For $v$ in $\R$, we now want to define a process 
\begin{equation}
(\overline{C}_{t}^{(1),v},\overline{B}_{t}^{(1),v})_{t\geq v-2b}\text{ having the same transition as }(C_{t}^{(1)},B_{t}^{(1)})\text{ and being stationary.}\label{eq:def-B-barre-(1)-v}
\end{equation}
 We set $(\overline{C}_{v-2b}^{(1),v},\overline{B}_{v-2b}^{(1),v})$
such that it has the law $\eta_{2}$. As we have given its transition,
the process $(\overline{C}_{t}^{(1),v},\overline{B}_{t}^{(1),v})_{t\geq v-2b}$
is well defined in law. In addition, we suppose that it is independent
of all the other processes. By Fact \ref{fact:Theorem-3.3-p.151},
the process $(\overline{C}_{t}^{(1),v},\overline{B}_{t}^{(1),v})_{t\geq v-2b}$
is stationary. 

We define the renewal times of $\overline{B}^{(1),v}$ by:
\[
\overline{S}_{1}^{(1),v}=\inf\{t\geq v-2b\,:\,\overline{B}_{t+}^{(1),v}\neq\overline{B}_{t-}^{(1),v}\}\,,
\]
\[
\text{and by recurrence,}\,\overline{S}_{k}^{(1),v}=\inf\{t>\overline{S}_{k-1}^{(1),v}\,:\,\overline{B}_{t+}^{(1),v}\neq\overline{B}_{t-}^{(1),v}\}\,.
\]
We also define, for all $t$, 
\[
\overline{N}_{t}^{(1),v}=\inf\{j\,:\,\overline{S}_{j}^{(1),v}>t\}\,.\text{}
\]
As will be seen later, the processes $\overline{B}^{(1),v}$ and $\overline{B}^{(2),v}$
are used to define asymptotic quantities (see, for example, Proposition
\ref{prop:conv-q-pair}) and we need them to be defined on an interval
$[v,+\infty)$with $v$ possibly in $\R^{-}$. The process $\overline{B}^{(2),v}$
is defined below (Section \ref{subsec:Two-stationary-processes}).

\subsection{Tagged fragments conditioned to \label{subsec:Tagged-fragments-conditioned}split
up ($(\widehat{B}^{(1),v},\widehat{B}^{(2),v})$).}

For $v$ in $[0,+\infty)$, we define a process $(\widehat{B}_{t}^{(1),v},\widehat{B}_{t}^{(2),v})_{t\geq0}$
such that
\begin{multline}
\widehat{B}^{(1),v}=B^{(1)}\text{ and, with }B^{(1)}\text{ fixed, }\widehat{B}^{(2),v}\text{ has the law of }B^{(2)}\text{ conditioned on }\\
\forall u\in\mathcal{U}\,,\,1\in A_{u}\Rightarrow[2\in A_{u}\Leftrightarrow-\log(\xi_{u})\leq v]\,,\label{eq:def-B-chapeau}
\end{multline}
which reads as follows~: the tag $2$ remains on the fragment bearing
the tag $1$ until the size of the fragment is smaller than $e^{-v}$.
We observe that, conditionally on $\widehat{B}_{v}^{(1),v}$, $\widehat{B}_{v}^{(2),v}$:
$(\widehat{B}_{v+\widehat{B}_{v}^{(1),v}+t}^{(1),v})_{t\geq0}$ and
$(\widehat{B}_{v+\widehat{B}_{v}^{(2),v}+t}^{(2),v})_{t\geq0}$ are
independent. We also define $\widehat{C}^{(1),v}=C^{(1)}$. There
is an algorithmic way to define $\widehat{B}^{(1),v}$ and $\widehat{B}^{(2),v}$,
which is illustrated in Figure \ref{fig:Processes-,-.}. Remember
$\widehat{B}^{(1),v}=B^{(1)}$and the definition of the mark $(\xi_{u},A_{u},I_{u})_{u\in\mathcal{U}}$
in Section \ref{subsec:Tagged-fragments}. 
\begin{figure}[h]
\begin{centering}
\includegraphics[scale=0.75]{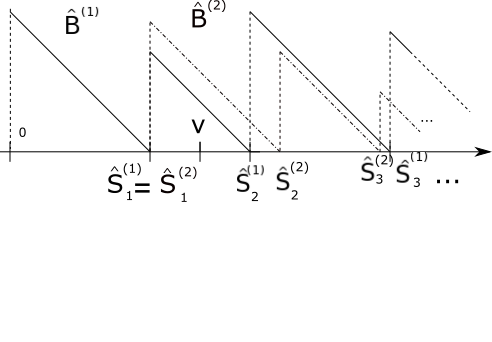}\caption{\label{fig:Processes-,-.}Processes $\widehat{B}^{(1),v}$, $\widehat{B}^{(2),v}$.}
\par\end{centering}
\end{figure}
We call $(\widehat{S}_{j}^{(i)})_{i=1,2;j\geq1}$ the renewal times
of these processes (as before, they can be defined as the times when
the right-hand side and left-hand side limits are not the same). If
$\widehat{S}_{j}^{(1)}\leq v$ then $\widehat{S}_{j}^{(2)}=\widehat{S}_{j}^{(1)}$.
If $k$ is such that $\widehat{S}_{k-1}^{(1)}\leq v$ and $\widehat{S}_{k}^{(1)}>v$,
we remember that 
\begin{equation}
\exp(\widehat{S}_{k}^{(1)}-\widehat{S}_{k-1}^{(1)})=\widetilde{\xi}_{ui}\label{eq:i-choose-i}
\end{equation}
 for some $u$ in $\mathcal{U}$ with $|u|=k-1$ and some $i$ in
$\N^{*}$ (because $\widehat{B}^{(1),v}=B^{(1)}$). We have points
$Y_{1},$$Y_{2}\in[0,1]$ such that $Y_{1},Y_{2}$ are in $I_{u}$
of length $\xi_{u}$%
. Conditionally on $\{Y_{1},Y_{2}\in I_{u}\}$, $Y_{1}$ and $Y_{2}$
are independent and uniformly distributed on $I_{u}$. The interval
$I_{ui}$, of length $\xi_{u}\widetilde{\xi}_{ui}$, is a sub-interval
of $I_{u}$ such that $Y_{1}\in I_{ui}$ (because of Equation (\ref{eq:i-choose-i})
above). Then, for $r\in\N^{*}\backslash\{i\}$, we want $Y_{2}$ to
be in $I_{ur}$ with probability $\widetilde{\xi}_{ur}/(1-\widetilde{\xi}_{ui})$
(because we want: $2\notin A_{ui}$). So we take 
\[
\widehat{S}_{k}^{(2)}=\widehat{S}_{k-1}^{(1)}-\log\widetilde{\xi}_{ur}
\]
 with probability $\widetilde{\xi}_{ur}/(1-\widetilde{\xi}_{ui})$
($r\in\N^{*}\backslash\{i\}$).
\begin{fact}
\label{fact:C-B}
\begin{enumerate}
\item The knowledge of the couple $(\widehat{S}_{N_{v}^{(1)}-1}^{(1)},\widehat{B}_{v}^{(1),v})$
is equivalent to the knowledge of the couple $(\widehat{C}_{v}^{(1),v},\widehat{B}_{v}^{(1),v})$.
\item The law of $B_{v}^{(1)}$ knowing $C_{v}^{(1)}$ is $\pi-C_{v}^{(1)}$
with $\pi$ conditioned to be bigger than $C_{v}^{(1)}$, we call
it $\eta_{1}(\dots|C_{v}^{(1)})$. As $\widehat{B}^{(1),v}=B^{(1)}$
and $\widehat{C}^{(1),v}=C^{(1)}$, we also have that the law of $\widehat{B}_{v}^{(1),v}$
knowing $\widehat{C}_{v}^{(1),v}$ is $\eta_{1}(\dots|\widehat{C}_{v}^{(1),v})$. 
\item The law of $\widehat{B}_{v}^{(2),v}$ knowing $(\widehat{C}_{v}^{(1),v},\widehat{B}_{v}^{(1),v})$
does not depend on $v$ and we denote it by $\eta'(\dots|\widehat{C}_{v}^{(1),v},\widehat{B}_{v}^{(1),v})$.
\end{enumerate}
\end{fact}

The subsequent waiting-time $\widehat{S}_{k+1}^{(1)}-\widehat{S}_{k}^{(1)}$,
$\widehat{S}_{k+1}^{(2)}-\widehat{S}_{k}^{(2)}$, \dots{} are chosen
independently of each other, each of them having the law $\pi$. For
$j$ equal to $1$ or $2$ and $t$ in $[0,+\infty)$, we define
\[
\widehat{N}_{t}^{(j)}=\inf\{i\,:\,\widehat{S}_{i}^{(j)}>t\}\,.
\]
We observe that for $t\geq2b$, $\widehat{N}_{t}^{(1)}$ is bigger
than $2$ (because of Assumption \ref{hyp:delta-step}).

\subsection{Two stationary processes after a split-up ($\overline{B}^{(1),v},\overline{B}^{(2),v}$)\label{subsec:Two-stationary-processes}.}

Let $k$ be an integer bigger than $2$ ($k\geq2$) and such that
\begin{equation}
k\times(b-a)\geq b\,.\label{eq:number-of-steps}
\end{equation}
 Now we state a small Lemma that will be useful below. Remember that
$(\overline{C}_{t}^{(1),v},\overline{B}_{t}^{(1),v})_{t\geq v-2b}$
is defined in Equation (\ref{eq:def-B-barre-(1)-v}). The process
$(\widehat{C}_{t}^{(1),v},\widehat{B}_{t}^{(1),v})_{t\geq0}$ is defined
in the previous Section. 
\begin{lem}
\label{lem:support}Let $v$ be in $\R.$ The variables $(\overline{C}_{v}^{(1),v},\overline{B}_{v}^{(1),v})$
and $(\widehat{C}_{kb}^{(1),kb},\widehat{B}_{kb}^{(1),kb})$ have
the same support (and it is $\mathcal{C}$, defined below Lemma \ref{lem:The-measure-}).
\end{lem}

\begin{proof}
The law $\eta_{2}$ is the law of $(\overline{C}_{0}^{(1)},\overline{B}_{0}^{(1)})$
($\eta_{2}$ is defined below Lemma \ref{lem:The-measure-}). As said
before, the support of $\eta_{2}$ is $\mathcal{C}$; and so (by stationarity)
the support of $(\overline{C}_{v}^{(1),v},\overline{B}_{v}^{(1),v})$
is $\mathcal{C}$.

Keep in mind that $\widehat{B}^{(1),v}=B^{(1)}$, $\widehat{C}^{(1),v}=C^{(1)}$.
By Assumption \ref{hyp:delta-step}, the support of $S_{k}^{(1)}$
is $[ka,kb]$ and the support of $S_{k+1}^{(1)}-S_{k}^{(1)}$ is $[a,b]$.
If $S_{k+1}^{(1)}>kb$ then $B_{kb}^{(1)}=S_{k+1}^{(1)}-S_{k}^{(1)}-(kb-S_{k}^{(1)})$
and $C_{kb}^{(1)}=kb-S_{k}^{(1)}$ (see Figure \ref{fig:B-et-C}).
\begin{figure}[h]
\centering{}\includegraphics[scale=0.7]{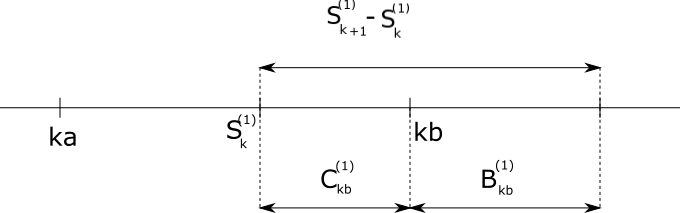}\caption{$B_{kb}^{(1)}$ and $C_{kb}^{(1)}$\label{fig:B-et-C}}
\end{figure}
The support of $S_{k}^{(1)}$is $[ka,kb]$ and $kb-ka\geq b$ (Equation
(\ref{eq:number-of-steps})), so, as $S_{k}^{(1)}$ and $S_{k+1}^{(1)}-S_{k}^{(1)}$
are independent, we get that the support of $(C_{kb}^{(1)},S_{k+1}^{(1)}-S_{k}^{(1)})$
includes $\{(u,w)\in[0;b]^{2}\,:\,w\geq\sup(a,u)\}$. And so, the
support of $(C_{kb}^{(1)},B_{kb}^{(1)})=(C_{kb}^{(1)},S_{k+1}^{(1)}-S_{k}^{(1)}-C_{kb}^{(1)})$
includes $\mathcal{C}$. As this support is included in $\mathcal{C}$,
we have proved the desired result.
\end{proof}
For $v$ in $\R$, we define a process $(\overline{B}_{t}^{(2),v})_{t\geq v}$.
We start by: 
\begin{multline}
\overline{B}_{v}^{(2),v}\text{ has the law }\eta'(\dots|\overline{C}_{v}^{(1),v},\overline{B}_{v}^{(1),v})\label{eq:def-B-barre-v-(2)}
\end{multline}
(remember $\eta'$ is defined in Fact \ref{fact:theta}). This conditioning
is correct because the law of $(\overline{C}_{v}^{(1),v},\overline{B}_{v}^{(1),v})$
is $\eta_{2}$ and its support is included in the support of the law
of $(\widehat{C}_{kb}^{(1),kb},\widehat{B}_{kb}^{(1),kb})$, which
is $\eta_{2}$ (see the Lemma above and below Equation (\ref{eq:def-B-barre-(1)-v})).%
{} We then let the process $(\overline{B}_{t}^{(1),v},\overline{B}_{t}^{(2),v})_{t\geq v}$
run its course as a Markov process having the same transition as $(\widehat{B}_{t-v+kb}^{(1),kb},\widehat{B}_{t-v+kb}^{(2),kb})_{t\geq v}$.
This means that, after time $v$ , $\overline{B}_{t}^{(1),v}$ and
$\overline{B}_{t}^{(2),v}$ decrease linearly (with a slope $-1$).
Until they reach $0$. When they reach $0$, each of these two processes
makes a jump of law $\pi$, independently of the other one. After
what, they decrease linearly \dots{} and so on.
\begin{fact}
\label{fact:By-construction,-the}The process $(\overline{B}_{t}^{(1),v},\overline{B}_{t}^{(2),v})_{t\geq v}$
is supposed independent from all the other processes (until now, we
have defined its law and said that that $\overline{B}^{(1),v}$ is
independent from all the other processes). 
\end{fact}

\section{\label{sec:Rate-of-convergence}Rate of convergence in the Key Renewal
Theorem}

We need the following regularity assumption. 

\begin{hypothesis}\label{hyp:queue-pi}

The probability $\pi(dx)$ is absolutely continuous with respect to
the Lebesgue measure (we will write $\pi(dx)=\pi(x)dx$). The density
function $x\mapsto\pi(x)$ is continuous on $(0;+\infty)$. 

\end{hypothesis}
\begin{fact}
\label{fact:theta}Let $\theta>1$ ($\theta$ is fixed in the rest
of the paper). The density $\pi$ satisfies 
\[
\limsup_{x\rightarrow+\infty}\exp(\theta x)\pi(x)<+\infty\,.
\]
\end{fact}

For $\varphi$ a nonnegative Borel-measurable function on $\R$, we
set $S(\varphi)$ to be the set of complex-valued measures $\rho$
(on the Borelian sets) such that $\int_{\R}\varphi(x)|\rho|(dx)<\infty$,
where $|\rho|$ stands for the total variation norm. If $\rho$ is
a finite complex-valued measure on the Borelian sets of $\R$, we
define $\mathcal{T}\rho$ to be the $\sigma$-finite measure with
the density 
\[
v(x)=\begin{cases}
\rho((x,+\infty)) & \text{ if }x\geq0\,,\\
-\rho((-\infty,x]) & \text{ if }x<0\,.
\end{cases}
\]
Let $F$ be the cumulative distribution function of $\pi$. 

We set $B_{t}=B_{t}^{(1)}$ (see Equation (\ref{eq:def-B}) for the
definition of $B^{(1)}$, $B^{(2)}$, \ldots ). By Theorem 3.3 p.
151 and Theorem 4.3 p. 156 of \cite{asmussen-2003}, we know that
$B_{t}$ converges in law to a random variable $B_{\infty}$(of law
$\eta$) and that $C_{t}$ converges in law to a random variable $C_{\infty}$
(of law $\eta$). The following Theorem is a consequence of \cite{sgibnev-2002},
Theorem 5.1, p. 2429. It shows there is actually a rate of convergence
for these convergences in law. 
\begin{thm}
\label{lem:sgibnev}Let $\epsilon'\in(0,\theta)$ . Let $M\in(0,+\infty)$.
Let 
\[
\varphi(x)=\begin{cases}
e^{(\theta-\epsilon')x} & \text{ if }x\geq0\,,\\
1 & \text{ if }x<0\,.
\end{cases}
\]
 %
If $Y$ is a random variable of law $\pi$ then
\begin{equation}
\sup_{\alpha\,:\,\Vert\alpha\Vert_{\infty}\leq M}\left|\E(\alpha(B_{t}))-\frac{1}{\mu}\int_{\R^{+}}\E(\alpha(Y-s)\1_{\{Y-s>0\}})ds\right|=o\left(\frac{1}{\varphi(t)}\right)\label{eq:sgibnev-B}
\end{equation}
as $t$ approaches $+\infty$ outside a set of Lebesgue measure zero
(the supremum is taken on $\alpha$ in the set of Borel-measurable
functions on $\R$), and 
\begin{equation}
\sup_{\alpha\,:\,\Vert\alpha\Vert_{\infty}\leq M}\left|\E(\alpha(C_{t}))-\frac{1}{\mu}\int_{\R^{+}}\E(\alpha(Y-s)\1_{\{Y-s>0\}})ds\right|=o\left(\frac{1}{\varphi(t)}\right)\label{eq:sgibnev-C}
\end{equation}
as $t$ approaches $+\infty$ outside a set of Lebesgue measure zero
(the supremum is taken on $\alpha$ in the set of Borel-measurable
functions on $\R$) .
\end{thm}

\begin{proof}
We write the proof of Equation (\ref{eq:sgibnev-B}) only. The proof
of Equation (\ref{eq:sgibnev-C}) is very similar. Let $*$ stands
for the convolution product. We define the renewal measure $U(dx)=\sum_{n=0}^{+\infty}\pi^{*n}(dx)$
(notations: $\pi^{*0}(dx)=\delta_{0}$, the Dirac mass at $0$, $\pi^{*n}=\pi*\pi*\dots*\pi$
($n$ times)).  We take i.i.d. variables $X,X_{1},X_{2}\dots$ of
law $\pi$. Let $f:\,\R\rightarrow\R$ be a measurable function such
that $\Vert f\Vert_{\infty}\leq M$. We have, for all $t\geq0$, 
\begin{eqnarray*}
\E(f(B_{t})) & = & \E\left(\sum_{n=0}^{+\infty}f(X_{1}+X_{2}+\dots+X_{n+1}-t)\1_{\{X_{1}+\dots+X_{n}\leq t<X_{1}+\dots+X_{n+1}\}}\right)\\
 & = & \int_{0}^{t}\E(f(s+X-t)\1_{\{s+X-t>0\}})U(ds)\,.
\end{eqnarray*}
We set 
\[
g(t)=\begin{cases}
\E(f(X-t)\1_{\{X-t>0\}}) & \text{ if }t\geq0\,,\\
0 & \text{ if }t<0\,.
\end{cases}
\]
We observe that $\Vert g\Vert_{\infty}\leq M$. %
{} We have, for all $t\geq0$,
\begin{eqnarray*}
\left|\E(f(X-t)\1_{\{X-t>0\}})\right| & \leq & \Vert f\Vert_{\infty}\p(X>t)\\
 & \leq & \Vert f\Vert_{\infty}e^{-(\theta-\frac{\epsilon'}{2})t}\E(e^{(\theta-\frac{\epsilon'}{2})X})\,.
\end{eqnarray*}
We have (by Fact \ref{fact:theta}): $\E(e^{(\theta-\frac{\epsilon'}{2})X})<\infty$.
The function $\varphi$ is submultiplicative and it is such that 
\[
\lim_{x\rightarrow-\infty}\frac{\log(\varphi(x))}{x}=0\leq\lim_{x\rightarrow+\infty}\frac{\log(\varphi(x))}{x}=\theta-\epsilon'\,.
\]
 The function $g$ is in $L^{1}(\R)$. The function $g.\varphi$ is
in $L^{\infty}(\R)$. We have $g(x)\varphi(x)\rightarrow0$ as $|x|\rightarrow\infty$.
We have
\[
\varphi(t)\int_{t}^{+\infty}|g(x)|dx\underset{t\rightarrow+\infty}{\longrightarrow}0\,,\,\varphi(t)\int_{-\infty}^{t}|g(x)|dx\underset{t\rightarrow-\infty}{\longrightarrow}0\,.
\]
We have $\mathcal{T}^{\circ2}(\pi)\in S(\varphi)$. 

Let us now take  a function $\alpha$ such that $\Vert\alpha\Vert_{\infty}\leq M$.
We set
\[
\widehat{\alpha}(t)=\begin{cases}
\E(\alpha(X-t)\1_{\{X-t\geq0\}}) & \text{ if }t\geq0\,,\\
0 & \text{ if }t<0\,.
\end{cases}
\]
Then we have $\Vert\widehat{\alpha}\Vert_{\infty}\leq M$ and (computing
as above for $f$)
\begin{eqnarray*}
\E(\alpha(B_{t})) & = & \widehat{\alpha}*U(t)
\end{eqnarray*}

In the case where $f$ is constant equal to $M$, we have $\Vert g\Vert_{\infty}=M$.
So, by \cite{sgibnev-2002}, Theorem 5.1 (applied to the case $f\equiv M$),
we have proved the desired result. 
\end{proof}
\begin{cor}
\label{cor:sgibnev}There exists a constant $\Gamma_{1}$ bigger than
$1$ such that: for any bounded measurable function $F$ on $\R$
such that $\eta(F)=0$,
\begin{equation}
|\E(F(B_{t}))|\leq\Vert F\Vert_{\infty}\times\frac{\Gamma_{1}}{\varphi(t)}\label{eq:cor-sgibnev-B}
\end{equation}
for $t$ outside a set of Lebesgue measure zero, and 
\begin{equation}
|\E(F(C_{t}))|\leq\Vert F\Vert_{\infty}\times\frac{\Gamma_{1}}{\varphi(t)}\label{eq:cor-sgibnev-C}
\end{equation}
for $t$ outside a set of Lebesgue measure zero.
\end{cor}

\begin{proof}
We write the proof of Equation (\ref{eq:cor-sgibnev-B}) only. The
proof of Equation (\ref{eq:cor-sgibnev-C}) is very similar. We take
$M=1$ in the above Theorem. Keep in mind  that $\eta$ is defined
in Equation (\ref{eq:def-eta}). By the above Theorem, there exists
a constant $\Gamma_{1}$ such that: for all measurable function $\alpha$
such that $\Vert\alpha\Vert_{\infty}\leq1$, 
\begin{equation}
\left|\E(\alpha(B_{t}))-\eta(\alpha)\right|\leq\frac{\Gamma_{1}}{\varphi(t)}\,\text{(for }t\text{ outside a set of Lebesque measure zero).}\label{eq:rate-01}
\end{equation}
Let us now take a bounded measurable $F$ such that $\eta(F)=0$.
By Equation (\ref{eq:rate-01}), we have (for $t$ outside a set of
Lebesgue measure zero)
\begin{eqnarray*}
\left|\E\left(\frac{F(B_{t})}{\Vert F\Vert_{\infty}}\right)-\eta\left(\frac{F}{\Vert F\Vert_{\infty}}\right)\right| & \leq & \frac{\Gamma_{1}}{\varphi(t)}\\
|\E(F(B_{t}))| & \leq & \Vert F\Vert_{\infty}\times\frac{\Gamma_{1}}{\varphi(t)}\,.
\end{eqnarray*}
\end{proof}

\section{\label{sec:Limits-of-symmetric}Limits of symmetric functionals }

\subsection{Notations\label{subsec:Notations-1}}

\begin{figure}[h]
\begin{centering}
\includegraphics[scale=0.75]{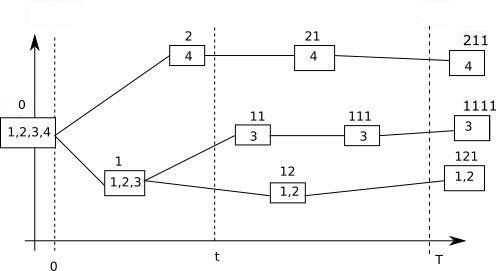}
\par\end{centering}
\caption{\label{fig:Tree-and-marks}Tree and marks}
\end{figure}

We fix $q\in\N^{*}$. We set $\mathcal{S}_{q}$ to be the symmetric
group of order $q$. A function $F\,:\,\R^{q}\rightarrow\R$ is symmetric
if
\[
\forall\sigma\in\mathcal{S}_{q}\,,\,\forall(x_{1},\dots,x_{q})\in\R^{q}\,,\,F(x_{\sigma(1)},x_{\sigma(2)},\dots,x_{\sigma(q)})=F(x_{1},x_{2},\dots,x_{q})\,.
\]
For $F:\R^{q}\rightarrow\R$, we define a symmetric version of $F$
by
\begin{equation}
F_{\text{sym}}(x_{1},\dots,x_{q})=\frac{1}{q!}\sum_{\sigma\in\mathcal{S}_{q}}F(x_{\sigma(1)},\dots,x_{\sigma(q)})\,,\,\text{for all }(x_{1},\dots,x_{q})\in\R^{q}\,.\label{eq:def-F-sym}
\end{equation}
We set $\mathcal{B}_{\text{sym }}(q)$ to be the set of bounded, measurable,
symmetric functions $F$ on $\R^{q}$, and we set $\mathcal{B}_{\text{sym }}^{0}(q)$
to be the $F$ of $\mathcal{B}_{\text{sym }}(q)$ such that 
\[
\int_{x_{1}}F(x_{1},x_{2},\dots,x_{q})\eta(dx_{1})=0\,,\,\forall(x_{2},\dots,x_{q})\in\R^{q-1}\,.
\]
Suppose that $k$ is in $[q]$ and $l\geq1$. %
For $t$ in $[0,T]$, we consider the following collections of nodes
of $\mathcal{U}$ (remember $T=-\log\epsilon$ and $\mathcal{U}$,
$\mathbf{m}(.)$ defined in Section \ref{subsec:Fragmentation-chains})~:
\[
\mathcal{T}_{1}=\{u\in\mathcal{U}\backslash\{0\}\,:\,A_{u}\neq\emptyset\,,\,\xi_{\boldsymbol{m}(u)}\geq\epsilon\}\cup\{0\}\,,
\]
\begin{equation}
S(t)=\{u\in\mathcal{T}_{1}\,:\,-\log(\xi_{\boldsymbol{m}(u)})\leq t\,,\,-\log(\xi_{u})>t\}=\mathcal{U}_{e^{-t}}\,,\label{eq:def_S(t)}
\end{equation}
\begin{equation}
L_{t}=\sum_{u\in S(t)\,:\,A_{u}\neq\emptyset}(\#A_{u}-1)\,.\label{eq:def-L_t}
\end{equation}
We set $\mathcal{L}_{1}$ to be the set of leaves in the tree $\mathcal{T}_{1}$.
For $t$ in $[0,T]$ and $i$ in $[q]$, there exists one and only
one $u$ in $S(t)$ such that $i\in A_{u}$. We call it $u\{t,i\}$.
Under Assumption \ref{hyp:delta-step}, there exists a constant bounding
the numbers vertices of $\mathcal{T}_{1}$ almost surely. Let us look
at an example in Figure \ref{fig:Tree-and-marks}. 
Here, we have a graphic representation of a realization of $\mathcal{T}_{1}$.
Each node $u$ of $\mathcal{T}_{1}$ is written above a rectangular
box in which we read $A_{u}$; the right side of the box has the coordinate
$-\log(\xi_{u})$ on the $X$-axis. For simplicity, the node $(1,1)$
is designated by $11$, the node $(1,2)$ is designated by $12$,
and so on. In this example: 
\begin{itemize}
\item $\mathcal{T}_{1}=\{(0),(1),(2),(1,1),(2,1),(1,2),(1,1,1),(2,1,1),(1,1,1,1),(1,2,1)\}$, 
\item $\mathcal{L}_{1}=\{(2,1,1),(1,1,1,1),(1,2,1)\}$, 
\item $A_{(1)}=\{1,2,3\}$, $A_{(1,2)}=\{1,2\}$, \dots{}
\item $S(t)=\{(1,2),(1,1),(2,1)\}$, 
\item $u\{t,1\}=(1,2)$, $u\{t,2\}=(1,2)$, $u\{t,3\}=(1,1)$, $u\{t,4\}=(2,1)$.
\end{itemize}
For $k$, $l$ in $\N$ and $t$ in $[0,T]$, we define the event
\[
C_{k,l}(t)=\{\sum_{u\in S(t)}\1_{\#A_{u}=1}=k\,,\,\sum_{u\in S(t)}(\#A_{u}-1)=l\}\,.
\]
For example, in Figure \ref{fig:Tree-and-marks}, we are in the event
$C_{2,1}(t)$. 

We define
\[
\mathcal{T}_{2}=\{u\in\mathcal{T}_{1}\backslash\{0\}\,:\,\#A_{\boldsymbol{m}(u)}\geq2\}\cup\{0\}\,,
\]
\[
m_{2}\,:\,u\in\mathcal{T}_{2}\mapsto(\xi_{u},\inf\{i,i\in A_{u}\})\,.
\]
For example, in Figure \ref{fig:Tree-and-marks}, $\mathcal{T}_{2}=\{(0),(1),(2),(1,1),(1,2),(1,2,1)\}$.
Let $\alpha$ be in $(0,1)$.
\begin{fact}
\label{fact:We-can-always-suppose}We can always suppose that $(1-\alpha)T>b$
because we are interested in $T$ going to infinity. So, in the following,
we suppose $(1-\alpha)T>b$.
\end{fact}

For any $t$, we can compute $\sum_{u\in S(t)}(\#A_{u}-1)$ if we
know $\sum_{u\in S(t)}\1_{\#A_{u}=1}$ and $\#S(t)$. As $T-\alpha T>b$,
any $u$ in $S(\alpha T)$ satisfies $\#A_{u}\geq2$ if and only if
$u$ is the mother of some $v$ in $\mathcal{T}_{2}$. So we deduce
that $C_{k,l}(\alpha T)$ is measurable with respect to $(\mathcal{T}_{2},m_{2})$.
We set, for all $u$ in $\mathcal{T}_{2}$, 
\begin{equation}
T_{u}=-\log(\xi_{u})\,.\label{eq:def-T_u}
\end{equation}
For any $i$ in $[q]$, $t\mapsto u\{t,i\}$ is piece-wise constant
and the the ordered sequence of its jump times is $S_{_{1}}^{(i)}<S_{2}^{(i)}<\dots$
(the $S_{\dots}^{(i)}$ are defined in Section \ref{subsec:Observation-scheme}).
We simply have that $1$, $e^{-S_{1}^{(i)}}$, $e^{-S_{2}^{i)}}$,
\dots{} are the successive sizes of the fragment supporting the tag
$i$. For example, in Figure \ref{fig:Tree-and-marks}, we have 
\begin{equation}
S_{1}^{(1)}=-\log(\xi_{1})\,,\,S_{2}^{(1)}=-\log(\xi_{(1,2)})\,,\,S_{3}^{(1)}=-\log(\xi_{(1,2,1)})\,,\dots\label{eq:def-T-(i)-j}
\end{equation}
 Let $\mathcal{L}_{2}$ be the set of leaves $u$ in the tree $\mathcal{T}_{2}$
such that the set $A_{u}$ has a single element $n_{u}$. For example,
in Figure \ref{fig:Tree-and-marks}, $\mathcal{L}_{2}=\{(2),(1,1)\}$.
We observe that $\#\mathcal{L}_{1}=q\Leftrightarrow\#\mathcal{L}_{2}=q$,
and thus
\begin{equation}
\{\mathcal{L}_{1}=q\}\in\sigma(\mathcal{L}_{2})\,.\label{eq:L1(q)_in_sigma...}
\end{equation}

We summarize the definition of $n_{u}$ in the following equation
\begin{equation}
\#A_{u}=1\Rightarrow A_{u}=\{n_{u}\}\,.\label{eq:def-n_u}
\end{equation}

For $q$ even ($q=2p$) and for all $t$ in $[0,T]$, we define the
events
\[
G_{t}=\{\forall i\in[p]\,,\,\exists u_{i}\in\mathcal{U}\,:\,\xi_{u_{i}}<e^{-t}\,,\,\xi_{\boldsymbol{m}(u_{i})}\geq e^{-t}\,,\,A_{u_{i}}=\{2i-1,2i\}\}\,,
\]

\[
\forall i\in[p]\,,\,G_{i,i+1}(t)=\{\exists u\in S(t)\,:\,\{2i-1,2i\}\subset A_{u}\}\,.
\]
We set, for all $t$ in $[0,T]$, 
\[
\mathcal{F}_{S(t)}=\sigma(S(t),(\xi_{u},A_{u})_{u\in S(t)})\,.
\]

\subsection{Intermediate results \label{subsec:Intermediate-results}}

The reader has to keep in mind that $T=-\log(\epsilon)$ (Equation
(\ref{eq:relation-T-epsilon})) and that $\delta$ is defined in Assumption
\emph{\ref{hyp:delta-step}. }The set $\Bsym(q)$ is defined in Section
\ref{subsec:Notations-1}.
\begin{lem}
\label{lem:conv-nouille-libre}We suppose that $F$ is in $\Bsym(q)$
and that $F$ is of the form $F=(f_{1}\otimes f_{2}\otimes\dots\otimes f_{q})_{\text{sym}}$,
with $f_{1}$, $f_{2}$, \dots{} , $f_{q}$ in $\Bsym(1)$ . Let $A$
be in $\sigma(\mathcal{L}_{2})$. For any $\alpha$ in $]0,1[$, $k$
in $[q]$ and $l$ in $\{0,1,\dots,(q-k-1)_{+}\}$, we have
\[
|\E(\1_{C_{k,l}(\alpha T)}\1_{A}F(B_{T}^{(1)},B_{T}^{(2)},\dots,B_{T}^{(q)}))|\leq\Vert F\Vert_{\infty}\Gamma_{1}^{q}C_{tree}(q)\left(\frac{1}{\delta}\right)^{q}\epsilon^{q/2}\,,
\]
(for a constant $C_{tree}(q)$ defined below in the proof and $\Gamma_{1}$
defined in Corollary \ref{cor:sgibnev}) and
\[
\]
\[
\epsilon^{-q/2}\E(\1_{C_{k,l}(\alpha T)}\1_{A}F(B_{T}^{(1)},B_{T}^{(2)},\dots,B_{T}^{(q)}))\underset{\epsilon\rightarrow0}{\longrightarrow}0\,.
\]
\end{lem}

\begin{proof}
{} Let $A$ be in $\sigma(\mathcal{L}_{2})$. %
We have

\begin{multline*}
\E(\1_{C_{k,l}(\alpha T)}\1_{A}F(B_{T}^{(1)},B_{T}^{(2)},\dots,B_{T}^{(q)}))\\
=\E(\1_{A}\sum_{f:\mathcal{T}_{2}\rightarrow\mathcal{P}([q])\text{ s.t. \dots}}\E(F(B_{T}^{(1)},B_{T}^{(2)},\dots,B_{T}^{(q)})\1_{A_{u}=f(u),\forall u\in\mathcal{T}_{2}}|\mathcal{L}_{2},\mathcal{T}_{2},m_{2}))\,\\
\text{(}\mathcal{P}\text{ defined in Section \ref{subsec:Notations})}
\end{multline*}
where we sum on the $f:\mathcal{T}_{2}\rightarrow\mathcal{P}([q])$
\uline{such that} 
\begin{equation}
\begin{cases}
f(u)=\sqcup_{v:\boldsymbol{m}(v)=u}f(v), & \text{ for all u in \ensuremath{\mathcal{T}_{2}},}\\
\sum_{u\in S(\alpha T)}\1_{\#f(u)=1}=k\text{ and }\sum_{u\in S(\alpha T)}(\#f(u)-1)=l
\end{cases}\label{eq:f-such-that}
\end{equation}
 We remind the reader that $\sqcup$ is defined in Section \ref{subsec:Notations}
(disjoint union), $\boldsymbol{m}(\dots)$ is defined in Section \ref{subsec:Fragmentation-chains}
(mother), $S(\dots)$ is defined in Equation (\ref{eq:def_S(t)}).
Here, we mean that we sum over the $f$ compatible with a description
of tagged fragments. 

If $u$ in $\mathcal{L}_{2}$ and if $T_{u}<T$, then, conditionally
on $\mathcal{T}_{2}$, $m_{2}$, $B_{T}^{(n_{u})}$ is independent
of all the other variables and has the same law as $B_{T-T_{u}}^{(1)}$
($T_{u}$ is defined in Equation (\ref{eq:def-T_u}), $n_{u}$ is
defined in Equation (\ref{eq:def-n_u})). Thus, using Theorem \ref{lem:sgibnev}
and Corollary \ref{cor:sgibnev}, we get, for any $\epsilon'\in(0,\theta-1)$,
$u\in\mathcal{L}_{2}$, 
\begin{multline*}
|\E(f_{n_{u}}(B_{T}^{(n_{u})})|\mathcal{L}_{2},\mathcal{T}_{2},m_{2})|\leq\Gamma_{1}\Vert f_{n_{u}}\Vert_{\infty}e^{-(\theta-\epsilon')(T-T_{u})_{+}}\,,\,\\
\text{for }T-T_{u}\notin Z_{0}\text{ where \ensuremath{Z_{0}} is of Lebesgue measure zero.}
\end{multline*}
Thus we get
\begin{multline*}
\frac{|\E(\1_{C_{k,l}(\alpha T)}\1_{A}F(B_{T}^{(1)},B_{T}^{(2)},\dots,B_{T}^{(q)}))|}{\Vert F\Vert_{\infty}\Gamma_{1}^{q}}\\
\text{(since }F\text{ is of the form }F=(f_{1}\otimes\dots\otimes f_{q})_{\text{sym}}\text{,}\\
\text{since, conditionally on }u\in\mathcal{L}_{2}\text{, }\text{the distribution of }T_{u}\text{ is absolutely continuous }\\
\text{with respect to the Lebesgue measure)}\\
\text{ }\\
\leq\E\left(\sum_{f:\mathcal{T}_{2}\rightarrow\mathcal{P}([q])\text{ s.t. \dots}}\left[\prod_{u\in\mathcal{L}_{2}}e^{-(\theta-\epsilon')(T-T_{u})_{+}}\times\1_{A}\E(\1_{A_{u}=f(u),\forall u\in\mathcal{T}_{2}}|\mathcal{L}_{2},\mathcal{T}_{2},m_{2})\right]\right)\\
\text{(because of Assumption \ref{hyp:delta-step} and because }\theta-\epsilon'>1\text{)}\\
\leq\E\left(\sum_{f:\mathcal{T}_{2}\rightarrow\mathcal{P}([q])\text{ s.t. \dots}}\left[\prod_{u\in\mathcal{L}_{2}}e^{-(T-T_{\boldsymbol{m}(u)})-\log(\delta)}\times\1_{A}\E(\1_{A_{u}=f(u),\forall u\in\mathcal{T}_{2}}|\mathcal{L}_{2},\mathcal{T}_{2},m_{2})\right]\right)\\
\text{(because of Equation (\ref{eq:proba-rester-dans-fragment}), see full proof in Section \ref{subsec:Detailed-proof-of})}\\
\leq\E\left(\sum_{f:\mathcal{T}_{2}\rightarrow\mathcal{P}([q])\text{ s.t. \dots}}\1_{A}\left[\prod_{u\in\mathcal{L}_{2}}e^{-(T-T_{\boldsymbol{m}(u)})-\log(\delta)}\times\prod_{u\in\mathcal{T}_{2}\backslash\{0\}}e^{-(\#f(u)-1)(T_{u}-T_{\boldsymbol{m}(u)})}\right]\right)\,.
\end{multline*}
For a fixed $\omega$ and a fixed $f$, we have
\[
\prod_{u\in\mathcal{L}_{2}}e^{-(T-T_{\boldsymbol{m}(u)})-\log(\delta)}\times\prod_{u\in\mathcal{T}_{2}\backslash\{0\}}e^{-(\#f(u)-1)(T_{u}-T_{\boldsymbol{m}(u)})}=\left(\frac{1}{\delta}\right)^{\#\mathcal{L}_{2}}\exp\left(-\int_{0}^{T}a(s)ds\right)\,,
\]
where, for all $s$, 
\begin{eqnarray*}
a(s) & = & \sum_{u\in\mathcal{L}_{2}\backslash\{0\}\,:\,T_{\boldsymbol{m}(u)}\leq s<T}\1_{\#f(u)=1}+\sum_{u\in\mathcal{T}_{2}\backslash\{0\}\,:\,T_{\boldsymbol{m}(u)}\leq s\leq T_{u}}(\#f(u)-1)\\
 &  & \text{(if }u\in\mathcal{T}_{2}\backslash\mathcal{L}_{2}\text{, }\1_{\#f(u)=1}=0\text{) }\\
 & = & \sum_{u\in\mathcal{T}_{2}\backslash\{0\}\,:\,T_{\boldsymbol{m}(u)}\leq s<T}\1_{\#f(u)=1}+\sum_{u\in\mathcal{T}_{2}\backslash\{0\}\,:\,T_{\boldsymbol{m}(u)}\leq s\leq T_{u}}(\#f(u)-1)\\
 &  & \text{(}S(.)\text{ defined in Equation (\ref{eq:def_S(t)}))}\\
 & \geq & \sum_{u\in S(s)}\1_{\#f(u)=1}+\sum_{u\in S(s)}(\#f(u)-1)\,.
\end{eqnarray*}
We observe that, under Equation (\ref{eq:f-such-that}): 
\[
a(t)\geq\left\lceil \frac{q}{2}\right\rceil \,,\,\forall t\,,
\]
\[
a(\alpha T)\geq k+l,
\]
and if $t$ is such that 
\[
\sum_{u\in S(t)}\1_{\#f(u)=1}=k'\,,\,\sum_{u\in S(t)}(\#f(u)-1)=l'\,,
\]
 for some integers $k'$, $l'$, then for all $s\geq t$,
\[
a(s)\geq k'+\left\lceil \frac{q-k'}{2}\right\rceil \,.
\]
We observe that, under Assumption \ref{hyp:delta-step}, there exists
a constant which bounds $\#\mathcal{T}_{2}$ almost surely (because
for all $u$ in $\mathcal{U}\backslash\{0\}$, $-\log(\xi_{u})+\log(\xi_{\mathbf{m}(u)})\geq a$)
and so there exists a constant $C_{tree}(q)$ which bounds $\#\{f:\mathcal{T}_{2}\rightarrow\mathcal{P}([q])\}$
almost surely. So, we have 
\begin{multline}
|\E(\1_{A}F(B_{T}^{(1)},B_{T}^{(2)},\dots,B_{T}^{(q)}))|\\
\leq\Vert F\Vert_{\infty}\Gamma_{1}^{q}\E\left(\sum_{f:\mathcal{T}_{2}\rightarrow\mathcal{P}([q])\text{ s.t. \dots}}\1_{A}\left(\frac{1}{\delta}\right)^{\#\mathcal{L}_{2}}e^{-\lceil q/2\rceil\alpha T}e^{-(k+\left\lceil \frac{q-k}{2}\right\rceil )(T-\alpha T)}\right)\\
\leq\Vert F\Vert_{\infty}\Gamma_{1}^{q}C_{tree}(q)\left(\frac{1}{\delta}\right)^{q}e^{-\lceil q/2\rceil\alpha T}e^{-(k+\left\lceil \frac{q-k}{2}\right\rceil )(1-\alpha)T}\,.\label{eq:borne-nouilles-libres}
\end{multline}
As $k\geq1$, then $k+\left\lceil \frac{q-k}{2}\right\rceil >\frac{q}{2}$,
and so we have proved the desired result (remember that $T=-\log\epsilon$). 
\end{proof}
\begin{rem}
\label{rem:About-Assumption} If we replaced Assumption \ref{hyp:delta-step}
by Assumption \ref{hyp:no-step}, we would have difficulties adapting
the above proof. In the second line of Equation (\ref{eq:borne-nouilles-libres})
above, the $1/\delta$ becomes $e^{T_{u}-T_{\boldsymbol{m}(u)}}$.
In addition, the tree $\mathcal{T}_{2}$ is not a.s. finite anymore.
So the expectation on the second line of (\ref{eq:borne-nouilles-libres})
could certainly be bounded, but for a high price (a lot more computations,
maybe assumptions on the tails of $\pi$, and so on). This is why
we stick to Assumption \ref{hyp:delta-step}. 
\end{rem}

Remember that $L_{t}$ ($t\geq0$) is defined in Equation (\ref{eq:def-L_t}).
\begin{lem}
\label{lem:conv-nouilles-liees}Let $k$ be an integer $\geq q/2$.
Let $\alpha\in[q/(2k),1]$. We have
\[
\p(L_{\alpha T}\geq k)\leq K_{1}(q)\epsilon^{q/2}\,,
\]
where $K_{1}(q)=\sum_{i\in[q]}\frac{q!}{(q-i)!}\times i^{q-i}$. 

Let $k$ be an integer $>q/2$. Let $\alpha\in(q/(2k),1)$. We have

\[
\epsilon^{-q/2}\p(L_{\alpha T}\geq k)\underset{\epsilon\rightarrow0}{\longrightarrow}0\,.
\]
(We remind the reader  that $T=-\log(\epsilon)$.) 
\begin{proof}
Let $k$ be an integer $\geq q/2$ and let $\alpha\in[q/(2k),1]$.
Remember that $S(\dots)$ is defined in Equation (\ref{eq:def_S(t)}).
Observe that: $\#S(\alpha T)=i$ if, and only if $L_{\alpha T}=q-i$
(see Equation (\ref{eq:def-L_t})). We decompose
\begin{eqnarray*}
\{L_{\alpha T}\geq k\} & = & \{L_{\alpha T}\in\{k,k+1,\dots,q-1\}\}\\
 & = & \cup_{i\in[q-k]}\{\#S(\alpha T)=i\}\\
 & = & \cup_{i\in[q-k]}\cup_{m:[i]\hookrightarrow[q]}(F(i,m)\cap\{\#S(\alpha T)=i\})\,,
\end{eqnarray*}
(remember that ``$\hookrightarrow$'' means we are summing on injections,
see Section \ref{subsec:Notations}) where
\[
F(i,m)=\{i_{1},i_{2}\in[i]\text{ with }i_{1}\ne i_{2}\Rightarrow\exists u_{1},u_{2}\in S(\alpha T),u_{1}\ne u_{2},m(i_{1})\in A_{u_{1}}\,,\,m(i_{2})\in A_{u_{2}}\}
\]
(to make the above equations easier to understand, observe that if
$\#S(\alpha T)=i$, we have for each $j\in[i]$, an index $m(j)$
in $A_{u}$ for some $u\in S(\alpha T)$, and we can choose $m$ such
that we are in the event $F(i,m)$). Suppose we are in the event $F(i,m)$.
For $u\in S(\alpha T)$ and for all $j$ in $[i]$ such that $m(j)\in A_{u}$,
we define (remember $|u|$ and $\mathbf{m}$ are defined in Section
\ref{subsec:Fragmentation-chains})
\[
T_{|u|}^{(j)}=-\log(\xi_{u})\,,\,T_{|u|-1}^{(j)}=-\log(\xi_{\boldsymbol{m}(u)})\,,\,\dots\,,\,T_{1}^{(j)}=-\log(\xi_{\boldsymbol{m}^{\circ(|u|-1)}(u)})\,,T_{0}^{(j)}=0\,,
\]
\[
l(j)=|u|\,,\,v(j)=u\,.
\]
We have
\begin{multline*}
\p(L_{\alpha T}\geq k)\leq\sum_{i\in[q-k]}\sum_{m:[i]\hookrightarrow[q]}\p(F(i,m)\cap\{\#S(\alpha T)=i\})\\
=\sum_{i\in[q-k]}\sum_{m:[i]\hookrightarrow[q]}\E(\1_{F(i,m)}\E(\1_{\#S(\alpha T)=i}|F(i,m),(T_{p}^{(j)})_{j\in[i],p\in[l(j)]},(v(j))_{j\in[i]}))\\
\text{(below, we sum over the partitions \ensuremath{\mathcal{B}} of }[q]\backslash m([i])\text{ into }i\text{ subsets \ensuremath{\mathcal{B}_{1},\mathcal{B}_{2},\dots,\mathcal{B}_{i}})}\\
=\sum_{i\in[q-k]}\sum_{m:[i]\hookrightarrow[q]}\sum_{\mathcal{B}}\E(\1_{F(i,m)}\E(\prod_{j\in[i]}\prod_{r\in\mathcal{B}_{j}}\1_{r\in A_{v(j)}}|F(i,m),(T_{p}^{(j)})_{j\in[i],p\in[l(j)]},(v(j))_{j\in[i]}))\\
\text{(as }Y_{1},\dots,Y_{q}\text{ defined in Section \ref{subsec:Tagged-fragments} are independant)}\\
=\sum_{i\in[q-k]}\sum_{m:[i]\hookrightarrow[q]}\sum_{\mathcal{B}}\E(\1_{F(i,m)}\prod_{j\in[i]}\prod_{r\in\mathcal{B}_{j}}\E(\1_{r\in A_{v(j)}}|F(i,m),(T_{p}^{(j)})_{j\in[i],p\in[l(j)]},(v(j))_{j\in[i]}))\\
\mbox{(because of Equation (\ref{eq:proba-rester-dans-fragment}) and Equation (\ref{eq:def-T_u}))}\\
=\sum_{i\in[q-k]}\sum_{m:[i]\hookrightarrow[q]}\sum_{\mathcal{B}}\E\left(\1_{F(i,m)}\prod_{j\in[i]}\prod_{r\in\mathcal{B}_{j}}\prod_{s=1}^{l(j)}\exp((-T_{s}^{(j)}+T_{s-1}^{(j)}))\right)\\
\text{(as }v(j)\in S(\alpha T)\text{)}\\
\leq\sum_{i\in[q-k]}\sum_{m:[i]\hookrightarrow[q]}\sum_{\mathcal{B}}\prod_{j\in[i]}\prod_{r\in\mathcal{B}_{j}}e^{-\alpha T}\\
=\sum_{i\in[q-k]}\sum_{m:[i]\hookrightarrow[q]}\sum_{\mathcal{B}}e^{-\alpha(q-i)T}\\
\leq\sum_{i\in[q-k]}\sum_{m:[i]\hookrightarrow[q]}\sum_{\mathcal{B}}e^{-k\alpha T}\leq e^{-k\alpha T}\times\sum_{i\in[q]}\frac{q!}{(q-i)!}i^{q-i}\,.
\end{multline*}
If we suppose that $k>q/2$ and $\alpha\in(q/(2k),1)$, then
\[
\exp\left(\frac{qT}{2}\right)\exp(-k\alpha T)\underset{T\rightarrow+\infty}{\longrightarrow}0\,.
\]
\end{proof}
\end{lem}

Immediate consequences of the two above lemmas are the following Corollaries.
\begin{cor}
\label{cor:lim-case-q-odd}If $q$ is odd and if $F\in\mathcal{B}_{\text{sym}}^{0}(q)$
is of the form $F=(f_{1}\otimes\dots\otimes f_{q})_{\text{sym}}$,
then%
\[
\]
\[
\epsilon^{-q/2}\E(F(B_{T}^{(1)},\dots,B_{T}^{(q)})\1_{\#\mathcal{L}_{1}=q})\underset{\epsilon\rightarrow0}{\longrightarrow}0\,.
\]
($\mathcal{B}_{\text{sym}}^{0}$ and $\mathcal{L}_{1}$ are defined
in Section \ref{subsec:Notations-1}.)
\end{cor}

\begin{proof}
We take $\alpha\in\left(\frac{q}{2}\left\lceil \frac{q}{2}\right\rceil ^{-1},1\right)$.
We observe that for $k$ in $[q]$, $t$ in $(0,T)$, 
\[
\sum_{u\in S(t)}\1_{\#A_{u}=1}=k\Rightarrow\sum_{u\in S(t)}(\#A_{u}-1)\in\{0,1,\dots,(q-k-1)_{+}\}
\]
and ($L_{t}$ defined in Equation (\ref{eq:def-L_t}))
\[
\sum_{u\in S(t)}\1_{\#A_{u}=1}=0\Rightarrow L_{t}\geq\lceil\frac{q}{2}\rceil\,.
\]
So we can decompose
\begin{multline}
\epsilon^{-q/2}\left|\E(F(B_{T}^{(1)},\dots,B_{T}^{(q)})\1_{\#\mathcal{L}_{1}=q})\right|\\
=|\epsilon^{-q/2}\sum_{k\in[q]}\sum_{l\in\{0,1,\dots,(q-k-1)_{+}\}}\E(\1_{C_{k,l}(\alpha T)}\1_{\#\mathcal{L}_{1}=q}F(B_{T}^{(1)},\dots,B_{T}^{(q)}))\\
+\epsilon^{-q/2}\E(\1_{L_{\alpha T}\geq\left\lceil q/2\right\rceil }\1_{\#\mathcal{L}_{1}=q}F(B_{T}^{(1)},\dots,B_{T}^{(q)}))|\\
\text{(by (\ref{eq:L1(q)_in_sigma...}) and Lemmas \ref{lem:conv-nouille-libre}, \ref{lem:conv-nouilles-liees}) }\underset{\epsilon\rightarrow0}{\longrightarrow}0\,\label{eq:dec-01}
\end{multline}
($\mathcal{L}_{1}$ and $\mathcal{L}_{2}$ defined in Section \ref{subsec:Notations-1}).
\end{proof}
\begin{cor}
\label{cor:maj-fonction-centree}Suppose $F\in\mathcal{B}_{sym}^{0}(q)$
is of the form $F=(f_{1}\otimes\dots\otimes f_{q})_{\text{sym}}$.
Let $A$ in $\sigma(\mathcal{L}_{2})$. Then
\[
|\E(F(B_{T}^{(1)},\dots,B_{T}^{(q)})\1_{A})|\leq\Vert F\Vert_{\infty}\epsilon^{q/2}\left\{ K_{1}(q)+\Gamma_{1}^{q}C_{\text{tree}}(q)\left(\frac{1}{\delta}\right)^{q}(q+1)^{2}\right\} \,
\]
\end{cor}

\begin{proof}
We get (as in Equation (\ref{eq:dec-01}) above)
\begin{multline*}
|\E(F(B_{T}^{(1)},\dots,B_{T}^{(q)})\1_{A})|=|\E(F(B_{T}^{(1)},\dots,B_{T}^{(q)})\1_{A}(\1_{L_{\alpha T}\geq q/2}+\sum_{k'\in[q]}\sum_{0\leq l\leq(q-k'-1)_{+}}\1_{C_{k',l}(\alpha T)}))|\\
\text{(from Lemmas \ref{lem:conv-nouille-libre}, \ref{lem:conv-nouilles-liees})}\\
\leq\Vert F\Vert_{\infty}\epsilon^{q/2}\left\{ K_{1}(q)+\Gamma_{1}^{q}C_{\text{tree}}(q)\left(\frac{1}{\delta}\right)^{q}\sum_{k'\in[q]}1+(q-k'-1)_{+}\right\} 
\end{multline*}
and $\sum_{k'\in[q]}1+(q-k'-1)_{+}\leq(q+1)^{2}$ (see Section \ref{subsec:Detailed-proof-of-3}
in the Appendix for a detailed proof).
\end{proof}
We now want to find the limit of $\epsilon^{-q/2}\E(\1_{L_{T}\leq q/2}\1_{\#\mathcal{L}_{1}=q}F(B_{T}^{(1)},\dots,B_{T}^{(q)}))$
when $\epsilon$ goes to $0$, for $q$ even. First we need a technical
lemma. 

For any $i$, the process $(B_{t}^{(i)})$ has a stationary law (see
Theorem 3.3 p. 151 of \cite{asmussen-2003}). Let $B_{\infty}$ be
a random variable having this stationary law $\eta$ (it has already
appeared in Section \ref{sec:Rate-of-convergence}). We can always
suppose that it is independent of all the other variables.
\begin{fact}
\label{fact:From-now-on,}From now on, when we have an $\alpha$ in
$(0,1)$, we suppose that $\alpha T-\log(\delta)<(T+\alpha T)/2$,
$(T+\alpha T)/2-\log(\delta)<T$ (this is true if $T$ is large enough).
(The constant $\delta$ is defined in Assumption \ref{hyp:delta-step}.)
\end{fact}

\begin{lem}
\label{lem:terme-b} Let $f_{1}$ , $f_{2}$ be in $\mathcal{B}_{sym}^{0}(1)$.
Let $\alpha$ belong to $(0,1)$. Let $\epsilon'$ belong to $(0,\theta-1)$
($\theta$ is defined in Fact \ref{fact:theta}). We have
\begin{equation}
\int_{-\infty}^{-\log(\delta)}e^{-v}|\E(f_{1}(\overline{B}_{0}^{(1),v})f_{2}(\overline{B}_{0}^{(2),v}))|dv<\infty\label{eq:terme-b-01}
\end{equation}
and (a.s.), for $T$ large enough, 
\begin{multline}
\left|e^{T-\alpha T-B_{\alpha T}^{(1)}}\E(f_{1}\otimes f_{2}(B_{T,}^{(1)}B_{T}^{(2)})\1_{G_{1,2}(T)^{\complement}}|\mathcal{F}_{S(\alpha T)},G_{\alpha T})\right.\\
\left.-\int_{-\infty}^{-\log(\delta)}e^{-v}\E(\1_{v\leq\overline{B}_{0}^{(1),v}}f_{1}(\overline{B}_{0}^{(1),v})f_{2}(\overline{B}_{0}^{(2),v}))dv\right|\\
\leq\Gamma_{2}\Vert f_{1}\Vert_{\infty}\Vert f_{2}\Vert_{\infty}\exp\left(-(T-\alpha T)\left(\frac{\theta-\epsilon'-1}{2}\right)\right)\,,\label{eq:terme-b-02}
\end{multline}
where 
\[
\Gamma_{2}=\frac{\Gamma_{1}^{2}}{\delta^{2+2(\theta-\epsilon')}(2(\theta-\epsilon')-1)}+\frac{\Gamma_{1}}{\delta^{\theta-\epsilon'}}+\frac{\Gamma_{1}^{2}}{\delta^{2(\theta-\epsilon')}(2(\theta-\epsilon')-1)}\,.
\]
(The processes $B^{(1)}$, $B^{(2)}$, $\overline{B}^{(1),v}$ , $\overline{B}^{(2),v}$
are defined in Sections \ref{subsec:Observation-scheme}, \ref{subsec:Stationary-age-process},
\ref{subsec:Two-stationary-processes}.)
\end{lem}

\begin{proof}
We have, for all $s$ in $[\alpha T+B_{\alpha T}^{(1)},T]$, (because
of Equation (\ref{eq:proba-rester-dans-fragment}) and Equation (\ref{eq:def-T_u})) 

\[
\p(u\{s,2\}=u\{s,1\}|\mathcal{F}_{S(\alpha T)},G_{\alpha T},(S_{j}^{(1)})_{j\geq1})=\exp(-(s+B_{s}^{(1)}-(\alpha T+B_{\alpha T}^{(1)}))
\]
(we remind the reader that $u\{s,1\}$, $G_{1,2}$, \dots{} are defined
in Section \ref{subsec:Notations-1}, below Equation (\ref{eq:def-L_t})).
Let us introduce the breaking time $\tau_{1,2}$ between $1$ and
$2$ as a random variable having the following property: conditionally
on $\mathcal{F}_{S(\alpha T)},G_{\alpha T},(S_{j}^{(1)})_{j\geq1}$,
$\tau_{1,2}$ has the density
\[
s\in\R\mapsto\1_{[\alpha T+B_{\alpha T}^{(1)},+\infty)}(s)e^{-(s-(\alpha T+B_{T}^{(1)}))}
\]
(this is a translation of an exponential law). We have the equalities:
$\alpha T+B_{\alpha T}^{(1)}=S_{j_{0}}^{(1)}$ for some $j_{0}$,
$T+B_{T}^{(1)}=S_{i_{0}}^{(1)}$ for some $i_{0}$. Here, we have
to make a comment on the definitions of Section \ref{subsec:Notations-1}.
In Figure \ref{fig:Tree-and-marks}, we have: $-\log(\xi_{(1,2)})=S_{2}^{(1)}$
(as in Equation (\ref{eq:def-T-(i)-j})), $S(S_{2}^{(1)})=\{(1,2,1),(1,1,1),(2,1)\}$,
$u\{-\log(\xi_{(1,2)}),1\}=\{(1,2,1)\}$. It is important to understand
this example before reading what follows. The breaking time $\tau_{1,2}$%
{} has the following interesting property (for all $k\geq j_{0}$)
\begin{multline*}
\p(u\{S_{k}^{(1)},2\}\neq u\{S_{k}^{(1)},1\}|\mathcal{F}_{S(\alpha T)},G_{\alpha T},(S_{j}^{(1)})_{j\geq1})\\
=\p(\tau_{1,2}\in[\alpha T+B_{\alpha T}^{(1)},S_{k}^{(1)}]|\mathcal{F}_{S(\alpha T)},G_{\alpha T},(S_{j}^{(1)})_{j\geq1})\,.
\end{multline*}
Just because we can, we impose, for all $k\geq j_{0}$, conditionally
on $\mathcal{F}_{S(\alpha T)},G_{\alpha T},(S_{j}^{(1)})_{j\geq1}$,
\[
\{u\{S_{k}^{(1)},2\}\neq u\{S_{k}^{(2)},1\}\}=\{\tau_{1,2}\in[\alpha T+B_{\alpha T}^{(1)},S_{k}^{(1)}]\}\,.
\]
Now, let $v$ be in $[\alpha T+B_{\alpha T}^{(1)},T+B_{T}^{(1)}]$.
We observe that, for all $v$ in $[\alpha T+B_{\alpha T}^{(1)},T+B_{T}^{(1)}]$,
\begin{multline*}
\E(f_{1}\otimes f_{2}(B_{T,}^{(1)}B_{T}^{(2)})\1_{G_{1,2}(T)^{\complement}}|\mathcal{F}_{S(\alpha T)},G_{\alpha T},(S_{j}^{(1)})_{j\geq1},\tau_{1,2}=v)\\
\text{(because of Equation (\ref{eq:def-B-chapeau}))}\\
=\E(f_{1}\otimes f_{2}(\widehat{B}_{T,}^{(1),v}\widehat{B}_{T}^{(2),v})|\mathcal{F}_{S(\alpha T)},G_{\alpha T},(S_{j}^{(1)})_{j\geq1})\,.
\end{multline*}

And so, 
\begin{multline*}
\E(f_{1}\otimes f_{2}(B_{T,}^{(1)}B_{T}^{(2)})\1_{G_{1,2}(T)^{\complement}}|\mathcal{F}_{S(\alpha T)},G_{\alpha T})\\
=\E\left(\E(f_{1}\otimes f_{2}(B_{T,}^{(1)}B_{T}^{(2)})\1_{G_{1,2}(T)^{\complement}}|\mathcal{F}_{S(\alpha T)},G_{\alpha T},(S_{j}^{(1)})_{j\geq1})|\mathcal{F}_{S(\alpha T)},G_{\alpha T}\right)\\
=\E(\E(\E(f_{1}\otimes f_{2}(B_{T,}^{(1)}B_{T}^{(2)})\1_{G_{1,2}(T)^{\complement}}|\mathcal{F}_{S(\alpha T)},G_{\alpha T},(S_{j}^{(1)})_{j\geq1},\tau_{1,2})\\
|\mathcal{F}_{S(\alpha T)},G_{\alpha T},(S_{j}^{(1)})_{j\geq1})|\mathcal{F}_{S(\alpha T)},G_{\alpha T})\\
\text{(keep in mind that }\widehat{B}^{(1),v}=B^{(1)}\text{ for all }v\text{)}\\
=\E(\E(\int_{\alpha T+B_{\alpha T}^{(1)}}^{T+B_{T}^{(1)}}e^{-(v-\alpha T-\widehat{B}_{\alpha T}^{(1),v})}\E(f_{1}\otimes f_{2}(B_{T,}^{(1)}B_{T}^{(2)})\1_{G_{1,2}(T)^{\complement}}|\mathcal{F}_{S(\alpha T)},G_{\alpha T},(S_{j}^{(1)})_{j\geq1},\tau_{1,2}=v)dv\\
|\mathcal{F}_{S(\alpha T)},G_{\alpha T},(S_{j}^{(1)})_{j\geq1})|\mathcal{F}_{S(\alpha T)},G_{\alpha T})\\
=\E(\E(\int_{\alpha T+B_{\alpha T}^{(1)}}^{T+B_{T}^{(1)}}e^{-(v-\alpha T-\widehat{B}_{\alpha T}^{(1),v})}\E(f_{1}(\widehat{B}_{T}^{(1),v})f_{2}(\widehat{B}_{T}^{(2),v})|\mathcal{F}_{S(\alpha T)},G_{\alpha T},(S_{j}^{(1)})_{j\geq1})dv\\
|\mathcal{F}_{S(\alpha T)},G_{\alpha T},(S_{j}^{(1)})_{j\geq1})|\mathcal{F}_{S(\alpha T)},G_{\alpha T})\\
=\E(\int_{\alpha T+B_{\alpha T}^{(1)}}^{T+B_{T}^{(1)}}e^{-(v-\alpha T-\widehat{B}_{\alpha T}^{(1),v})}f_{1}(\widehat{B}_{T}^{(1),v})f_{2}(\widehat{B}_{T}^{(2),v})dv|\mathcal{F}_{S(\alpha T)},G_{\alpha T})
\end{multline*}
Let us split the above integral into two parts and multiply them by
$e^{T-\alpha T-B_{\alpha T}^{(1)}}$. We have (this is the first part)
\begin{multline}
\left|e^{T-\alpha T-B_{\alpha T}^{(1)}}\E\left(\int_{\alpha T+B_{\alpha T}^{(1)}}^{(T+\alpha T)/2}e^{-(v-\alpha T-\widehat{B}_{\alpha T}^{(1),v})}f_{1}(\widehat{B}_{T}^{(1),v})f_{2}(\widehat{B}_{T}^{(2),v})dv|\mathcal{F}_{S(\alpha T)},G_{\alpha T}\right)\right|\\
=e^{T-\alpha T-B_{\alpha T}^{(1)}}|\E(\int_{\alpha T+B_{\alpha T}^{(1)}}^{(T+\alpha T)/2}e^{-(v-\alpha T-\widehat{B}_{\alpha T}^{(1),v})}\\
\times\E(f_{1}(\widehat{B}_{T}^{(1),v})f_{2}(\widehat{B}_{T}^{(2),v})|\widehat{B}_{v}^{(1),v},\widehat{B}_{v}^{(2),v},\mathcal{F}_{S(\alpha T)},G_{\alpha T})dv\\
|\mathcal{F}_{S(\alpha T)},G_{\alpha T})|\\
\text{(using the fact that }\widehat{B}_{T}^{(1),v}\text{ and }\widehat{B}_{T}^{(2),v}\text{ are independant }\\
\text{conditionally on }\{\widehat{B}_{v}^{(1),v},\widehat{B}_{v}^{(2),v},\mathcal{F}_{S(\alpha T)},G_{\alpha T}\}\text{ if }T\geq v-\log(\delta)\text{, }\\
\text{we get, by Theorem \ref{lem:sgibnev} and Corollary \ref{cor:sgibnev} and Fact \ref{fact:From-now-on,})}\\
\leq e^{T-\alpha T-B_{\alpha T}^{(1)}}\\
\times\E(\int_{\alpha T+B_{\alpha T}^{(1)}}^{(T+\alpha T)/2}e^{-(v-\alpha T-\widehat{B}_{\alpha T}^{(1)})}(\Gamma_{1}\Vert f_{1}\Vert_{\infty}e^{-(\theta-\epsilon')(T-v-\widehat{B}_{v}^{(1),v})_{+}}\times\Gamma_{1}\Vert f_{2}\Vert_{\infty}e^{-(\theta-\epsilon')(T-v-\widehat{B}_{v}^{(2),v})_{+}})dv\\
|\mathcal{F}_{S(\alpha T)},G_{\alpha T})\\
\text{(using Assumption \ref{hyp:delta-step})}\\
\leq\Gamma_{1}^{2}\Vert f_{1}\Vert_{\infty}\Vert f_{2}\Vert_{\infty}e^{T-\alpha T-\log(\delta)}\int_{\alpha T}^{(T+\alpha T)/2}e^{-(v-\alpha T+\log(\delta))}e^{-2(\theta-\epsilon')(T-v+\log(\delta))}dv\\
=\frac{\Gamma_{1}^{2}\Vert f_{1}\Vert_{\infty}\Vert f_{2}\Vert_{\infty}}{\delta^{2+2(\theta-\epsilon')}}e^{T-2(\theta-\epsilon')T}\left[\frac{e^{(2(\theta-\epsilon')-1)v}}{2(\theta-\epsilon')-1}\right]_{\alpha T}^{(T+\alpha T)/2}\\
\leq\frac{\Gamma_{1}^{2}\Vert f_{1}\Vert_{\infty}\Vert f_{2}\Vert_{\infty}}{\delta^{2+2(\theta-\epsilon')}}\frac{\exp\left(-(2(\theta-\epsilon')-1)T+(2(\theta-\epsilon')-1)\frac{(T+\alpha T)}{2}\right)}{2(\theta-\epsilon')-1}\\
=\frac{\Gamma_{1}^{2}\Vert f_{1}\Vert_{\infty}\Vert f_{2}\Vert_{\infty}}{\delta^{2+2(\theta-\epsilon')}}\frac{\exp\left(-(2(\theta-\epsilon')-1)\left(\frac{T-\alpha T}{2}\right)\right)}{2(\theta-\epsilon')-1}\,.\label{eq:bout-01}
\end{multline}
We have (this is the second part, minus some other term)
\begin{multline}
\left|\begin{array}{c}
\underbrace{e^{T-\alpha T-B_{\alpha T}^{(1)}}\E\left(\int_{(T+\alpha T)/2}^{T+B_{T}^{(1)}}e^{-(v-\alpha T-B_{\alpha T}^{(1)})}f_{1}(\widehat{B}_{T}^{(1),v})f_{2}(\widehat{B}_{T}^{(2),v})dv|\mathcal{F}_{S(\alpha T)},G_{\alpha T}\right)}\\
\text{second part}
\end{array}\right.\\
\left.-\begin{array}{c}
\underbrace{\int_{(T+\alpha T)/2}^{T-\log(\delta)}e^{-(v-T)}\E(\1_{v\leq T+\overline{B}_{T}^{(1),v}}f_{1}(\overline{B}_{T}^{(1),v})f_{2}(\overline{B}_{T}^{(2),v}))dv}\\
(\heartsuit)
\end{array}\right|\\
=\left|e^{T-\alpha T-B_{\alpha T}^{(1)}}\E(\int_{(T+\alpha T)/2}^{T-\log(\delta)}e^{-(v-\alpha T-B_{\alpha T}^{(1)})}\1_{v\leq T+B_{T}^{(1)}}f_{1}(\widehat{B}_{T}^{(1),v})f_{2}(\widehat{B}_{T}^{(2),v})dv|\mathcal{F}_{S(\alpha T)},G_{\alpha T})\right.\\
\left.-e^{T-\alpha T-B_{\alpha T}^{(1)}}\E(\int_{(T+\alpha T)/2}^{T-\log(\delta)}e^{-(v-\alpha T-B_{\alpha T}^{(1)})}\1_{v\leq T+\overline{B}_{T}^{(1),v}}f_{1}(\overline{B}_{T}^{(1),v})f_{2}(\overline{B}_{T}^{(2),v})dv)\right|\\
=e^{T-\alpha T-B_{\alpha T}^{(1)}}\left|\int_{(T+\alpha T)/2}^{T-\log(\delta)}e^{-(v-\alpha T-B_{\alpha T}^{(1)})}\E(\E(\1_{v\leq T+B_{T}^{(1)}}f_{1}(\widehat{B}_{T}^{(1),v})f_{2}(\widehat{B}_{T}^{(2),v})\right.\\
|\widehat{C}_{v}^{(1),v},\mathcal{F}_{S(\alpha T)},G_{\alpha T})|\mathcal{F}_{S(\alpha T)},G_{\alpha T})dv\\
-\int_{(T+\alpha T)/2}^{T-\log(\delta)}e^{-(v-\alpha T-B_{\alpha T}^{(1)})}\E(\E(\1_{v\leq T+\overline{B}_{T}^{(1),v}}f_{1}(\overline{B}_{T}^{(1),v})f_{2}(\overline{B}_{T}^{(2),v})\left.|\overline{C}_{v}^{(1),v}))dv\right|\label{eq:abs-a-borner-01}
\end{multline}
We observe that, for all $v$ in $[(T+\alpha T)/2,T-\log(\delta)]$,
once $\widehat{C}_{v}^{(1),v}$ is fixed, we can make a simulation
of $\widehat{B}_{T}^{(1),v}=B_{T}^{(1)}$, $\widehat{B}_{T}^{(2),v}$
(these processes are independent of $\mathcal{F}_{S(\alpha T)},G_{\alpha T}$
conditionally on $\widehat{C}_{v}^{(1),v}$). Indeed, we draw $\widehat{B}_{v}^{(1),v}$
conditionally on $\widehat{C}_{v}^{(1),v}$ (with law $\eta_{1}(\dots|\widehat{C}_{v}^{(1),v})$
defined in Fact \ref{fact:C-B}), then we draw $\widehat{B}_{v}^{(2),v}$conditionally
on $\widehat{B}_{v}^{(1),v}$ and $\widehat{C}_{v}^{(1),v}$ (with
law $\eta'(\dots|\widehat{B}_{v}^{(1),v},\widehat{C}_{v}^{(1),v})$,
see Fact \ref{fact:C-B}). Then, $(\widehat{B}_{t}^{(1),v})_{t\geq v}$,
$(\widehat{B}_{t}^{(2),v})_{t\geq v}$ run their courses as independent
Markov processes, until we get $\widehat{B}_{T}^{(1),v}$, $\widehat{B}_{T}^{(2),v}$. 

In the same way (for all $v$ in $[(T+\alpha T)/2,T-\log(\delta)]$),
we observe that the process $(\overline{C}^{(1),v},\overline{B}^{(1),v})$
starts at time $v-2b$ and has the same transition as $(C^{(1)},B^{(1)})$
(see Equation (\ref{eq:def-B-barre-(1)-v})). By Assumption \ref{hyp:A},
there exists the following time: $S=\sup\{t\,:\,v-b\leq t\leq v\,,\,\overline{C}_{t}^{(1),v}=0\}$.
We then have $v-S=\overline{C}_{v}^{(1),v}$. When $\overline{C}_{v}^{(1),v}$
is fixed, this entails that $\overline{B}_{v}^{(1),v}$ has the law
$\eta_{1}(\dots|\overline{C}_{v}^{(1),v})$. We have $\overline{B}_{v}^{(2),v}$
of law $\eta'(\dots|\overline{C}_{v}^{(1),v},\overline{B}_{v}^{(1),v})$
(by Equation (\ref{eq:def-B-barre-v-(2)})). As said before, we then
let the process $(\overline{B}_{t}^{(1),v},\overline{B}_{t}^{(2),v})_{t\geq v}$
run its course as a Markov process having the same transition as $(\widehat{B}_{t-v+kb}^{(1),kb},\widehat{B}_{t-v+kb}^{(2),kb})_{t\geq v}$
until we get $\overline{B}_{T}^{(1),v}$, $\overline{B}_{T}^{(2),v}$. 

So we get that (for all $v$ in $[(T+\alpha T)/2,T-\log(\delta)]$)
\[
\E(\1_{v\leq T+B_{T}^{(1)}}f_{1}(\widehat{B}_{T}^{(1),v})f_{2}(\widehat{B}_{T}^{(2),v})|\widehat{C}_{v}^{(1),v},\mathcal{F}_{S(\alpha T)},G_{\alpha T})=\Psi(\widehat{C}_{v}^{(1),v})\,,
\]
\[
\E(\1_{v\leq T+\overline{B}_{T}^{(1),v}}f_{1}(\overline{B}_{T}^{(1),v})f_{2}(\overline{B}_{T}^{(2),v})|\overline{C}_{v}^{(1),v})=\Psi(\overline{C}_{v}^{(1),v})\overset{\text{law}}{=}\Psi(C_{\infty})\,,
\]
for some function $\Psi$, the \uline{same} on both lines, such
that $\Vert\Psi\Vert_{\infty}\leq\Vert f_{1}\Vert_{\infty}\Vert f_{2}\Vert_{\infty}$($C_{\infty}$
defined in Section \ref{sec:Rate-of-convergence}). So, by Theorem
\ref{lem:sgibnev} and Corollary \ref{cor:sgibnev} applied on the
time interval $[\alpha T+B_{\alpha T}^{(1)},v]$, the quantity in
Equation (\ref{eq:abs-a-borner-01}) can be bounded by (remember that
$\widehat{C}^{(1),v}=C^{(1)}$, Section \ref{subsec:Tagged-fragments-conditioned})
\[
e^{T-\alpha T-B_{\alpha T}^{(1)}}\int_{(T+\alpha T)/2}^{T-\log(\delta)}e^{-(v-\alpha T-B_{\alpha T}^{(1)})}\Gamma_{1}\Vert f_{1}\Vert_{\infty}\Vert f_{2}\Vert_{\infty}e^{-(\theta-\epsilon')(v-\alpha T-B_{\alpha T}^{(1)})}dv
\]
(coming from Corollary \ref{cor:sgibnev} there is an integral over
a set of Lebesgue measure zero in the above bound, but this term vanishes).
The above bound can in turn be bounded by:
\begin{multline}
\Gamma_{1}\Vert f_{1}\Vert_{\infty}\Vert f_{2}\Vert_{\infty}\delta^{-(\theta-\epsilon')}e^{T}\int_{(T+\alpha T)/2}^{T-\log(\delta)}e^{(\theta-\epsilon')\alpha T}e^{-(\theta-\epsilon'+1)v}dv\\
\text{(as }\theta-\epsilon'+1>1\text{) }\\
\leq\Gamma_{1}\Vert f_{1}\Vert_{\infty}\Vert f_{2}\Vert_{\infty}\delta^{-(\theta-\epsilon')}e^{T+\alpha T(\theta-\epsilon')}\exp(-(\theta-\epsilon'+1)(\frac{T+\alpha T}{2}))\\
=\Gamma_{1}\Vert f_{1}\Vert_{\infty}\Vert f_{2}\Vert_{\infty}\delta^{-(\theta-\epsilon')}\exp\left(-(\theta-\epsilon'-1)\left(\frac{T-\alpha T}{2}\right)\right)\,.\label{eq:bout-02}
\end{multline}
We have
\begin{multline}
\int_{\frac{T+\alpha T}{2}}^{T-\log(\delta)}e^{-(v-T)}\E(\1_{v\leq T+\overline{B}_{T}^{(1),v}}f_{1}(\overline{B}_{T}^{(1),v})f_{2}(\overline{B}_{T}^{(2),v}))dv\\
\text{(as }(\overline{B}_{T}^{(1),v},\overline{B}_{T}^{(2),v})\text{ and }(\overline{B}_{0}^{(1),v-T},\overline{B}_{0}^{(2),v-T})\text{ have same law)}\\
=\int_{\frac{T+\alpha T}{2}}^{T-\log(\delta)}e^{-(v-T)}\E(\1_{v-T\leq\overline{B}_{0}^{(1),v-T}}f_{1}(\overline{B}_{0}^{(1),v-T})f_{2}(\overline{B}_{0}^{(2),v-T}))dv\\
\text{(change of variable }v'=v-T\text{)}\\
=\E\left(\int_{-\left(\frac{T-\alpha T}{2}\right)}^{-\log(\delta)}e^{-v'}\1_{v'\leq\overline{B}_{0}^{(1),v'}}f_{1}(\overline{B}_{0}^{(1),v'})f_{2}(\overline{B}_{0}^{(1),v'})dv'\right)\label{eq:bout-03}
\end{multline}
and
\begin{multline}
\int_{-\infty}^{-\frac{(T-\alpha T)}{2}}e^{-v}|\E(f_{1}(\overline{B}_{0}^{(1),v})f_{2}(\overline{B}_{0}^{(2),v}))|dv\\
\text{(since }\overline{B}_{0}^{(1),v}\text{ and }\overline{B}_{0}^{(2),v}\text{ are independant conditionnaly on }\overline{B}_{v}^{(1),v}\text{, }\overline{B}_{v}^{(2),v}\\
\text{if }v-\log(\delta)\leq0\text{)}\\
\text{(using Theorem \ref{lem:sgibnev} and Corollary \ref{cor:sgibnev}) }\\
\leq\int_{-\infty}^{-\frac{(T-\alpha T)}{2}}e^{-v}\Gamma_{1}^{2}\Vert f_{1}\Vert_{\infty}\Vert f_{2}\Vert_{\infty}\E(e^{-(\theta-\epsilon')(-v-\overline{B}_{v}^{(1),v})_{+}}e^{-(\theta-\epsilon')(-v-\overline{B}_{v}^{(2),v})_{+}})dv\\
\text{(again, coming from Corollary \ref{cor:sgibnev} there is an integral }\\
\text{over a set of Lebesgue measure zero in the above bound, but this term vanishes)}\\
\leq\int_{-\infty}^{-\frac{(T-\alpha T)}{2}}e^{-v}\Gamma_{1}^{2}\Vert f_{1}\Vert_{\infty}\Vert f_{2}\Vert_{\infty}e^{-2(\theta-\epsilon')(-v+\log(\delta))}dv\\
=\frac{\Gamma_{1}^{2}\Vert f_{1}\Vert_{\infty}\Vert f_{2}\Vert_{\infty}}{\delta^{2(\theta-\epsilon')}}\frac{\exp\left(-(2(\theta-\epsilon')-1)\frac{(T-\alpha T)}{2}\right)}{2(\theta-\epsilon')-1}\,.\label{eq:bout-04}
\end{multline}
Equations (\ref{eq:bout-03}) and (\ref{eq:bout-04}) give us Equation
(\ref{eq:terme-b-01}). Equations (\ref{eq:bout-01}), (\ref{eq:bout-02}),
(\ref{eq:bout-03}) and (\ref{eq:bout-04}) give us the desired result
(see Figure \ref{fig:Puzzle} to understand the puzzle). 
\begin{figure}[h]
\begin{centering}
\begin{tabular}{ccc}
 &  & Term $e^{T-\alpha T-B_{\alpha T}^{(1)}}\E(f_{1}\otimes f_{2}(B_{T,}^{(1)}B_{T}^{(2)})\1_{G_{1,2}(T)^{\complement}}|\mathcal{F}_{S(\alpha T)},G_{\alpha T})$\tabularnewline
 & $\swarrow$ & $\downarrow$\tabularnewline
First part (in Equation (\ref{eq:bout-01})) &  & Second part (in Equation (\ref{eq:abs-a-borner-01}))\tabularnewline
$\downarrow$ &  & $\downarrow$\tabularnewline
Small &  & Close to some term ($\heartsuit$) (see Equations (\ref{eq:abs-a-borner-01}),
(\ref{eq:bout-02}))\tabularnewline
 &  & $\downarrow$\tabularnewline
 &  & ($\heartsuit$) close to $-\int_{-\infty}^{-\log(\delta)}e^{-v}\E(\1_{v\leq\overline{B}_{0}^{(1),v}}f_{1}(\overline{B}_{0}^{(1),v})f_{2}(\overline{B}_{0}^{(2),v}))dv$\tabularnewline
 &  & (see Equations (\ref{eq:bout-03}), (\ref{eq:bout-04}))\tabularnewline
\end{tabular}
\par\end{centering}
\caption{Puzzle\label{fig:Puzzle}}

\end{figure}
\end{proof}
\begin{lem}
\label{lem:calcul-exact-limite}Let $k$ in $\{0,1,2,\dots,p\}$.
We suppose $q$ is even and $q=2p$. Let $\alpha\in(q/(q+2),1)$.
We suppose $F=f_{1}\otimes f_{2}\otimes\dots\otimes f_{q}$, with
$f_{1}$, \ldots{} , $f_{q}$ in $\mathcal{B}_{sym}^{0}(1)$. We then
have~:
\begin{multline}
\epsilon^{-q/2}\E(F(B_{T,}^{(1)}\dots,B_{T}^{(q)})\1_{G_{\alpha T}}\1_{\#\mathcal{L}_{1}=q})\\
\underset{\epsilon\rightarrow0}{\longrightarrow}\prod_{i=1}^{p}\int_{-\infty}^{-\log(\delta)}e^{-v}\E(\1_{v\leq\overline{B}_{0}^{(1),v}}f_{2i-1}(\overline{B}_{0}^{(1),v})f_{2i}(\overline{B}_{0}^{(2),v}))dv\,.\label{eq:conv-01}
\end{multline}
(Remember that $T=-\log\epsilon$.)
\end{lem}

\begin{proof}
By Fact \ref{fact:We-can-always-suppose}, we have $T>\alpha T-\log(\delta)$.
We have (remember the definitions just before Section \ref{subsec:Intermediate-results})
\[
G_{\alpha T}\cap\{\#\mathcal{L}_{1}=q\}=G_{\alpha T}\cap\underset{1\leq i\leq p}{\bigcap}G_{2i-1,2i}(T)^{\complement}\,.
\]
We have (remember $T=-\log(\epsilon)$)
\begin{multline}
\epsilon^{-q/2}\E(F(B_{T,}^{(1)}\dots,B_{T}^{(q)})\1_{G_{\alpha T}}\1_{\#\mathcal{L}_{1}=q})\\
=e^{pT}\E\left(\1_{G_{\alpha T}}\E\left(\left.\prod_{i=1}^{p}f_{2i-1}\otimes f_{2i}(B_{T}^{(2i-1)},B_{T}^{(2i)})\1_{G_{2i-1,2i}(T)^{\complement}}\right|\mathcal{F}_{S(\alpha T)},G_{\alpha T}\right)\right)\\
\text{(as }(B_{T}^{(1)},B_{T}^{(2)},\1_{G_{1,2}(T)})\text{, }(B_{T}^{(3)},B_{T}^{(4)},\1_{G_{3,4}(T)})\text{, \dots\ are independant conditionally on }\mathcal{F}_{S(\alpha T)},G_{\alpha T}\\
\text{due to Fact \ref{fact:From-now-on,})}\\
=\E\left(\1_{G_{\alpha T}}\prod_{i=1}^{p}e^{T}\E\left(\left.f_{2i-1}\otimes f_{2i}(B_{T}^{(2i-1)},B_{T}^{(2i)})\1_{G_{2i-1,2i}(T)^{\complement}}\right|\mathcal{F}_{S(\alpha T)},G_{2i-1,2i}(\alpha T)\right)\right)\\
\text{(by Lemma \ref{lem:terme-b} and as \ensuremath{(B^{(1)},\dots,B^{(q)})} is exchangeable)}\\
=\E\left(\1_{G_{\alpha T}}\prod_{i=1}^{p}e^{\alpha T+B_{\alpha T}^{(2i-1)}}\left(\int_{-\infty}^{-\log(\delta)}e^{-v}\E(\1_{v\leq\overline{B}_{0}^{(1),v}}f_{2i-1}(\overline{B}_{0}^{(1),v})f_{2i}(\overline{B}_{0}^{(2),v}))dv+R_{2i-1,2i}\right)\right)\,,\label{eq:conv-01-a}
\end{multline}
with (a.s.) 
\begin{equation}
|R_{2i-1,2i}|\leq\Gamma_{2}\Vert f_{2i-1}\Vert_{\infty}\Vert f_{2i}\Vert_{\infty}e^{-(T-\alpha T)\frac{(\theta-\epsilon'-1)}{2}}\,.\label{eq:conv-01-b}
\end{equation}
We introduce the events (for $t\in[0,T]$)($u\{\dots\}$ defined below
Equation (\ref{eq:def-L_t}))
\[
O_{t}=\left\{ \#\{u\{t,2i-1\},1\leq i\leq p\}=p\right\} \,,
\]
and the tribes (for $i$ in $[q]$, $t\in[0,T]$)
\[
\mathcal{F}_{t,i}=\sigma(u\{t,i\},\xi_{u\{t,i\}})\,.
\]
As $G_{\alpha T}=O_{\alpha T}\cap\bigcap_{1\leq i\leq p}\{u\{\alpha T,2i-1\}=u\{\alpha T,2i\}\}$,
we have~:
\begin{multline}
\E(\1_{G_{\alpha T}}\prod_{i=1}^{p}e^{B_{\alpha T}^{(2i-1)}+\alpha T})=\E(\1_{O_{\alpha T}}\prod_{i=1}^{p}e^{B_{\alpha T}^{(2i-1)}+\alpha T}\E(\prod_{i=1}^{p}\1_{u\{\alpha T,2i-1\}=u\{\alpha T,2i\}}|\vee_{1\leq i\leq p}\mathcal{F}_{\alpha T,2i-1}))\\
\text{(by Proposition \ref{prop:(reformulation-of-Proposition} and Equation (\ref{eq:proba-rester-dans-fragment}))}=\E(\1_{O_{\alpha T}})\,.\label{eq:1_O}
\end{multline}
We then observe that
\[
O_{\alpha T}^{\complement}=\cup_{i\in[p]}\cup_{j\in[p],j\neq i}\{u\{\alpha T,2i-1\}=u\{\alpha T,2j-1\}\}\,,
\]
and, for $i\neq j$,
\begin{eqnarray*}
\p(u\{\alpha T,2i-1\}=u\{\alpha T,2j-1\}) & = & \E(\E(\1_{u\{\alpha T,2i-1\}=u\{\alpha T,2j-1\}}|\mathcal{F}_{\alpha T,2i-1}))\\
\text{(by Proposition \ref{prop:(reformulation-of-Proposition} and Equation (\ref{eq:proba-rester-dans-fragment}))} & = & \E(e^{-\alpha T-B_{\alpha T}^{(2i-1)}})\\
\mbox{(because of Assumption (\ref{hyp:delta-step}))} & \leq & e^{-\alpha T-\log(\delta)}\,.
\end{eqnarray*}
So
\[
\p(O_{\alpha T})\underset{\epsilon\rightarrow0}{\longrightarrow}1\,.
\]
This gives us enough material to finish the proof of Equation (\ref{eq:conv-01}).

Indeed, starting from Equation (\ref{eq:conv-01-a}), we have
\begin{multline*}
\epsilon^{-q/2}\E(F(B_{T,}^{(1)}\dots,B_{T}^{(q)})\1_{G_{\alpha T}}\1_{\#\mathcal{L}_{1}=q})\\
=\E\left(\1_{G_{\alpha T}}\prod_{i=1}^{p}e^{\alpha T+B_{\alpha T}^{(2i-1)}}\right)\times\prod_{i=1}^{p}\left(\int_{-\infty}^{-\log(\delta)}e^{-v}\E(\1_{v\leq\overline{B}_{0}^{(1),v}}f_{2i-1}(\overline{B}_{0}^{(1),v})f_{2i}(\overline{B}_{0}^{(2),v}))dv\right)\\
+\E\left(\1_{G_{\alpha T}}\prod_{i=1}^{p}e^{\alpha T+B_{\alpha T}^{(2i-1)}}R_{2i-1,2i}\right)=:(1)+(2)\,.
\end{multline*}
By Equation  (\ref{eq:1_O}), 
\begin{multline*}
(1)=\p(O_{\alpha t})\prod_{i=1}^{p}\left(\int_{-\infty}^{-\log(\delta)}e^{-v}\E(\1_{v\leq\overline{B}_{0}^{(1),v}}f_{2i-1}(\overline{B}_{0}^{(1),v})f_{2i}(\overline{B}_{0}^{(2),v}))dv\right)\\
\underset{\epsilon\rightarrow0}{\longrightarrow}\prod_{i=1}^{p}\left(\int_{-\infty}^{-\log(\delta)}e^{-v}\E(\1_{v\leq\overline{B}_{0}^{(1),v}}f_{2i-1}(\overline{B}_{0}^{(1),v})f_{2i}(\overline{B}_{0}^{(2),v}))dv\right)\,.
\end{multline*}
And, by Equation (\ref{eq:conv-01-b}), 
\[
(2)\leq\p(O_{\alpha t})\prod_{i=1}^{p}\left(\Gamma_{2}\Vert f_{2i-1}\Vert_{\infty}\Vert f_{2i}\Vert_{\infty}e^{-(T-\alpha T)\frac{(\theta-\epsilon'-1)}{2}}\right)\underset{\epsilon\rightarrow0}{\longrightarrow}0\,.
\]
\end{proof}

\subsection{Convergence result}

For $f$ and $g$ bounded measurable functions, we set
\begin{equation}
V(f,g)=\int_{-\infty}^{-\log(\delta)}e^{-v}\E(\1_{v\leq\overline{B}_{0}^{(1),v}}f(\overline{B}_{0}^{(1),v})g(\overline{B}_{0}^{(2),v}))dv\,.\label{eq:def-V}
\end{equation}
For $q$ even, we set $\mathcal{I}_{q}$ to be the set of partitions
of $[q]$ into subsets of cardinality $2$. We have
\begin{equation}
\#\mathcal{I}_{q}=\frac{q!}{\left(q/2\right)!2^{q/2}}\,.\label{eq:card-I-q}
\end{equation}
 For $I$ in $\mathcal{I}_{q}$ and $t$ in $[0,T]$, we introduce
\[
G_{t,I}=\{\forall\{i,j\}\in I\,,\,\exists u\in\mathcal{U}\text{ such that }\xi_{u}<e^{-t}\,,\,\xi_{\boldsymbol{m}(u)}\geq e^{-t}\,,\,A_{u}=\{i,j\}\}\,.
\]
For $t$ in $[0,T]$, we define
\[
\mathcal{P}_{t}=\cup_{I\in\mathcal{I}_{q}}G_{t,I}\,.
\]
The above event can be understood as ``at time $t$, the dots are
paired on different fragments''. As before, the reader has to keep
in mind that $T=-\log(\epsilon)$ (Equation (\ref{eq:relation-T-epsilon})).
\begin{prop}
\label{prop:conv-q-pair}Let $q$ be in $\N^{*}$. Let $F=(f_{1}\otimes\dots\otimes f_{q})_{\text{sym }}$
with $f_{1}$, \dots , $f_{q}$ in $\Bsym(1)$ ($(\dots)_{\text{sym}}$
defined in Equation (\ref{eq:def-F-sym})). If $q$ is even ($q=2p$)
then
\begin{equation}
\epsilon^{q/2}\E(F(B_{T}^{(1)},\dots,B_{T}^{(q)})\1_{\#\mathcal{L}_{1}=q})\underset{\epsilon\rightarrow0}{\longrightarrow}\sum_{I\in\mathcal{I}_{q}}\prod_{\{a,b\}\in I}V(f_{a},f_{b})\,.\label{eq:conv-q-even}
\end{equation}
\end{prop}

\begin{proof}
Let $\alpha$ be in $(q/(q+2),1)$. We have
\begin{eqnarray*}
\epsilon^{-q/2}\E(F(B_{T}^{(1)},\dots,B_{T}^{(q)})\1_{\#\mathcal{L}_{1}=q}) & = & \epsilon^{-q/2}\E(F(B_{T}^{(1)},\dots,B_{T}^{(q)})\1_{\#\mathcal{L}_{1}=q}(\1_{\mathcal{P}_{\alpha T}}+\1_{\mathcal{P}_{\alpha T}^{\complement}}))\,.
\end{eqnarray*}
Remember that the events of the form $C_{k,l}(t)$, $L_{t}$ are defined
in Section \ref{subsec:Notations-1}. The set $\mathcal{P}_{\alpha T}^{\complement}$
is a disjoint union of sets of the form $C_{k,l}(\alpha T)$ (with
$k\geq1$) and $\{L_{\alpha T}>q/2\}$ (this can be understood heuristically
by: ``if the dots are not paired on fragments than some of them are
alone on their fragment, or none of them is alone on a fragment and
some are a group of at least three on a fragment''). As said before,
the event $\{\#\mathcal{L}_{1}=q\}$ is measurable with respect to
$\mathcal{L}_{2}$ (see Equation (\ref{eq:L1(q)_in_sigma...})). So,
by Lemma \ref{lem:conv-nouille-libre} and Lemma \ref{lem:conv-nouilles-liees},
we have that
\[
\lim_{\epsilon\rightarrow0}\epsilon^{-q/2}\E(F(B_{T}^{(1)},\dots,B_{T}^{(q)})\1_{\#\mathcal{L}_{1}=q}\1_{\mathcal{P}_{\alpha T}^{\complement}})=0\,.
\]
 We compute~:
\begin{multline*}
\epsilon^{-q/2}\E(F(B_{T}^{(1)},\dots,B_{T}^{(q)})\1_{\#\mathcal{L}_{1}=q}\1_{\mathcal{P}_{\alpha T}})=\epsilon^{-q/2}\E(F(B_{T}^{(1)},\dots,B_{T}^{(q)})\1_{\#\mathcal{L}_{1}=q}\sum_{I_{q}\in\mathcal{I}_{q}}\1_{G_{\alpha T,I_{q}}})\\
\text{(as }F\text{ is symmetric and }(B_{T}^{(1)},\dots,B_{T}^{(q)})\text{ is exchangeable)}\\
=\frac{q!}{2^{q/2}\left(\frac{q}{2}\right)!}\epsilon^{-q/2}\E(F(B_{T}^{(1)},\dots,B_{T}^{(q)})\1_{\#\mathcal{L}_{1}=q}\1_{G_{\alpha T}})\\
=\frac{q!\epsilon^{-q/2}}{2^{q/2}\left(\frac{q}{2}\right)!}\frac{1}{q!}\sum_{\sigma\in\mathcal{S}_{q}}\E((f_{\sigma(1)}\otimes\dots\otimes f_{\sigma(q)})(B_{T}^{(1)},\dots,B_{T}^{(q)})\1_{\#\mathcal{L}_{1}=q}\1_{G_{\alpha T}})\\
\text{(by Lemma \ref{lem:calcul-exact-limite}) }\\
\underset{\epsilon\rightarrow0}{\longrightarrow}\frac{1}{2^{q/2}\left(\frac{q}{2}\right)!}\sum_{\sigma\in\mathcal{S}_{q}}\prod_{i=1}^{p}V(f_{\sigma(2i-1)},f_{\sigma(2i)})=\sum_{I\in\mathcal{I}_{q}}\prod_{\{a,b\}\in I}V(f_{a},f_{b})\,.
\end{multline*}
\begin{eqnarray*}
\end{eqnarray*}
\end{proof}

\section{\label{sec:Results}Results}

We are interested in the probability measure $\gamma_{T}$ defined
by its action on bounded measurable functions $F\,:\,[0,1]\rightarrow\R$
by
\[
\gamma_{T}(F)=\sum_{u\in\mathcal{U}_{\epsilon}}\xi{}_{u}F\left(\frac{\xi{}_{u}}{\epsilon}\right)\,.
\]
We define, for all $q$ in $\N^{*}$, $F$ from $[0,1]^{q}$ to $\R$
,
\[
\gamma_{T}^{\otimes q}(F)=\sum_{a\,:\,[q]\rightarrow\mathcal{U_{\epsilon}}}\xi{}_{a(1)}\dots\xi{}_{a(q)}F\left(\frac{\xi_{a(1)}}{\epsilon},\dots,\frac{\xi_{a(q)}}{\epsilon}\right)\,,
\]

\[
\gamma_{T}^{\odot q}(F)=\sum_{a\,:\,[q]\hookrightarrow\mathcal{U_{\epsilon}}}\xi_{a(1)}\dots\xi_{a(q)}F\left(\frac{\xi_{a(1)}}{\epsilon},\dots,\frac{\xi_{a(q)}}{\epsilon}\right)\,,
\]
where the last sum is taken over all the injective applications $a$
from $[q]$ to $\Ue$. We set 
\[
\Phi(F)\,:\,(y_{1},\dots,y_{q})\in\R^{+}\mapsto F(e^{-y_{1}},\dots,e^{-y_{q}})\,,
\]

The law $\gamma^{\otimes q}$ is the law of $q$ fragments picked
in $\mathcal{U}_{\epsilon}$ with replacement. For each fragment,
the probability to be picked is its size. The measure $\gamma^{\odot q}$
is not a law: $\gamma^{\odot q}(F)$ is an expectation over $q$ fragments
picked in $\mathcal{U}_{\epsilon}$ with replacement (for each fragment,
the probability to be picked is its size), in this expectation, we
multiply the integrand by zero if two fragments are the same (and
by one otherwise). The definition of Section \ref{subsec:Tagged-fragments}
says that we can define the tagged fragment by painting colored dots
on the stick $[0,1]$ ($q$ dots of different colors, these are the
$Y_{1}$, \dots , $Y_{q}$) and then by looking on which fragments
of $\mathcal{U}_{\epsilon}$ we have these dots. So, we get (remember
$T=-\log\epsilon$)
\begin{equation}
\E(\gamma_{T}^{\otimes q}(F))=\E(\Phi(F)(B_{T}^{(1)},\dots,B_{T}^{(q)}))\,,\label{eq:lien-otimes-B}
\end{equation}

\begin{equation}
\E(\gamma_{T}^{\odot q}(F))=\E(\Phi(F)(B_{T}^{(1)},\dots,B_{T}^{(q)})\1_{\#\mathcal{L}_{1}}=q)\,.\label{eq:lien-odot-B}
\end{equation}
We define, for all bounded continuous $f\,:\,\R^{+}\rightarrow\R$,
\begin{equation}
\gamma_{\infty}(f)=\eta(\Phi(f))\,.\label{eq:def-gamma-infty}
\end{equation}

\begin{prop}
[Law of large numbers] \label{prop:convergence-ps}We remind the
reader that we have Fact \ref{fact:theta} and that we are under Hypothesis
\ref{hyp:A}, \ref{hyp:conservative}, \ref{hyp:delta-step}, \ref{hyp:queue-pi}.
Let f be a continuous function from $[0,1]$ to $\R$. We have:
\[
\gamma_{T}(f)\underset{T\rightarrow+\infty}{\overset{\text{a.s.}}{\longrightarrow}}\gamma_{\infty}(f)\,.
\]
(Remember $T=-\log\epsilon$.)
\end{prop}

\begin{proof}
We take a bounded measurable function $f\,:\,[0,1]\rightarrow\R$.
We define $\overline{f}=f-\eta(\Phi(f))$. We take an integer $q\geq2$.
We introduce the notation~:
\[
\forall g\,:\,\R^{+}\rightarrow\R\,,\,\forall(x_{1},\dots,x_{q})\in\R^{q}\,,\,g^{\otimes q}(x_{1},\dots,x_{q})=g(x_{1})g(x_{2})\dots g(x_{q})\,.
\]
 We have
\begin{eqnarray*}
\E((\gamma_{T}(f)-\eta(\Phi(f)))^{q}) & = & \E((\gamma_{T}(\overline{f}))^{q})\\
 & = & \E(\gamma_{t}^{\otimes q}(\overline{f}^{\otimes q}))\\
 &  & \text{(as }(B^{(1)},\dots,B^{(q)})\text{ is exchangeable)}\\
 & = & \E(\gamma_{t}^{\otimes q}((\overline{f}^{\otimes q})_{\text{sym}}))\\
\text{(by Corollary \ref{cor:maj-fonction-centree})} & \leq & \Vert\overline{f}\Vert_{\infty}^{q}\epsilon^{q/2}\left\{ K_{1}(q)+\Gamma_{1}^{q}C_{\text{tree}}(q)\left(\frac{1}{\delta}\right)^{q}(q+1)^{2}\right\} \,.
\end{eqnarray*}
We now take sequences $(T_{n}=\log(n))_{n\geq1}$, $(\epsilon_{n}=1/n)_{n\geq1}$.
We then have, for all $n$ and for all $\iota>0$,
\[
\p([\gamma_{T_{n}}(f)-\eta(\Phi(f))]^{4}\geq\iota)\leq\frac{\Vert\overline{f}\Vert_{\infty}^{4}}{\iota n^{2}}\left\{ K_{1}(4)+\Gamma_{1}^{4}C_{\text{tree}}(4)\left(\frac{1}{\delta}\right)^{4}\times25\right\} \,.
\]
So, by Borell-Cantelli's Lemma,
\begin{equation}
\gamma_{T_{n}}(f)\underset{n\rightarrow+\infty}{\underset{\longrightarrow}{\text{a.s.}}}\eta(\Phi(f))\,.\label{eq:conv-discrete}
\end{equation}

We now have to work a little bit more to get to the result. Let $n$
be in $\N^{*}$. We can decompose ($\mathcal{U}_{\epsilon}$ defined
in Section \ref{subsec:Observation-scheme}, $\sqcup$ stands for
``disjoint union'' and is defined in Section \ref{subsec:Notations})
\[
\mathcal{U}_{\epsilon_{n}}=\mathcal{U}_{\epsilon_{n}}^{(1)}\sqcup\mathcal{U}_{\epsilon_{n}}^{(2)}\text{ where }\mathcal{U}_{\epsilon_{n}}^{(1)}=\mathcal{U}_{\epsilon_{n}}\cap\mathcal{U}_{\epsilon_{n+1}}=\mathcal{U}_{\epsilon_{n+1}}\,,\,\mathcal{U}_{\epsilon_{n}}^{(2)}=\mathcal{U}_{\epsilon_{n}}\backslash\mathcal{U}_{\epsilon_{n+1}}\,.
\]
For $u$ in $\mathcal{U}_{\epsilon_{n}}\backslash\mathcal{U}_{\epsilon_{n+1}}$,
we set $\boldsymbol{d}(u)=\{v\,:\,u=\mathbf{\mathbf{m}}(v)\}$ ($\mathbf{m}$
is defined in Section \ref{subsec:Fragmentation-chains}) and we observe
that, for all $u$ ($T_{u}$ defined in Equation (\ref{eq:def-T_u}))
\begin{equation}
\sum_{v\in\boldsymbol{d}(u)}\xi_{v}=\xi_{u}\,.\label{eq:d-subset-a}
\end{equation}
 We can then write
\begin{equation}
\sum_{u\in\mathcal{U}_{\epsilon_{n}}}\xi_{u}f\left(\frac{\xi_{u}}{\epsilon_{n}}\right)=\sum_{u\in\mathcal{U}_{\epsilon_{n}}^{(1)}}\xi_{u}f(n\xi_{u})+\sum_{u\in\mathcal{U}_{\epsilon_{n}}^{(2)}}\xi_{u}f(n\xi_{u})\,.\label{eq:dec-U-epsilon-01}
\end{equation}
There exists $n_{1}$ such that, for $n$ bigger than $n_{1}$, $e^{-a}<\epsilon_{n+1}/\epsilon_{n}$
(remember Assumption \ref{hyp:delta-step}). We suppose $n\geq n_{1}$,
we then have, for all $u$ in $\mathcal{U}_{\epsilon_{n}}^{(2)}$,
$\epsilon_{n}>\xi_{u}\geq\epsilon_{n+1}$ and, for any $v$ in $\mathbf{d}(u)$,
$\xi_{v}\leq\epsilon_{n}e^{-a}$, $\xi_{v}<\epsilon_{n+1}$. So we
get 
\begin{equation}
\sum_{u\in\mathcal{U}_{\epsilon_{n+1}}}\xi_{u}f\left(\frac{\xi_{u}}{\epsilon_{n+1}}\right)=\sum_{u\in\mathcal{U}_{\epsilon_{n}}^{(1)}}\xi_{u}f((n+1)\xi_{u})+\sum_{u\in\mathcal{U}_{\epsilon_{n}}^{(2)}}\sum_{v\in\boldsymbol{d}(u)}\xi_{v}f((n+1)\xi_{v})\,.\label{eq:dec-U-epsilon-02}
\end{equation}
 So we have, for $n\geq n_{1}$, 
\begin{multline}
\left|\sum_{u\in\mathcal{U}_{\epsilon_{n}}^{(2)}}\sum_{v\in\boldsymbol{d}(u)}\xi_{v}f((n+1)\xi_{v})-\sum_{u\in\mathcal{U}_{\epsilon_{n}}^{(2)}}\xi_{u}f(n\xi_{u})\right|\\
\leq|\gamma_{T_{n+1}}(f)-\gamma_{T_{n}}(f)|+\left|\sum_{u\in\mathcal{U}_{\epsilon_{n}}^{(1)}}\xi_{u}f((n+1)\xi_{u})-\sum_{u\in\mathcal{U}_{\epsilon_{n}}^{(1)}}\xi_{u}f(n\xi_{u})\right|\,.\label{eq:terme-a}
\end{multline}
If we take $f=\Id$, the terms in the equation above can be bounded:
\begin{multline}
\left|\sum_{u\in\mathcal{U}_{\epsilon_{n}}^{(2)}}\sum_{v\in\boldsymbol{d}(u)}\xi_{v}f((n+1)\xi_{v})-\sum_{u\in\mathcal{U}_{\epsilon_{n}}^{(2)}}\xi_{u}f(n\xi_{u})\right|\\
\geq\left|\sum_{u\in\mathcal{U}_{\epsilon_{n}}^{(2)}}\left(\xi_{u}f(n\xi_{u})-\sum_{v\in\boldsymbol{d}(u)}\xi_{v}f(n\xi_{v})\right)\right|-\left|\sum_{u\in\mathcal{U}_{\epsilon_{n}}^{(2)}}\sum_{v\in\boldsymbol{d}(u)}(\xi_{v}f(n\xi_{v})-\xi_{v}f((n+1)\xi_{v}))\right|\\
\text{(by Assumption \ref{hyp:delta-step})}\\
\geq\sum_{u\in\mathcal{U}_{\epsilon_{n}}^{(2)}}\left(\xi_{u}f(n\xi_{u})-\sum_{v\in\boldsymbol{d}(u)}\xi_{v}f(n\xi_{u})e^{-a}\right)-\left|\sum_{u\in\mathcal{U}_{\epsilon_{n}}^{(2)}}\sum_{v\in\boldsymbol{d}(u)}(\xi_{v}f(n\xi_{v})-\xi_{v}f((n+1)\xi_{v}))\right|\\
\text{(by Equation (\ref{eq:d-subset-a})) }\geq\sum_{u\in\mathcal{U}_{\epsilon_{n}}^{(2)}}\xi_{u}(1-e^{-a})\frac{n}{n+1}-\sum_{u\in\mathcal{U}_{\epsilon_{n}}^{(2)}}\sum_{v\in\boldsymbol{d}(u)}\xi_{v}\frac{1}{n+1}\,,\label{eq:terme-b}
\end{multline}
\begin{multline}
|\gamma_{T_{n+1}}(f)-\gamma_{T_{n}}(f)|+\left|\sum_{u\in\mathcal{U}_{\epsilon_{n}}^{(1)}}\xi_{u}f((n+1)\xi_{u})-\sum_{u\in\mathcal{U}_{\epsilon_{n}}^{(1)}}\xi_{u}f(n\xi_{u})\right|\\
\leq|\gamma_{T_{n+1}}(f)-\gamma_{T_{n}}(f)|+\sum_{u\in\mathcal{U}_{\epsilon_{n}}^{(1)}}\xi_{u}\frac{1}{n}\,.\label{eq:terme-c}
\end{multline}
\uline{Let \mbox{$\iota>0$}}. We fix $\omega$ in $\Omega$. By
Equation (\ref{eq:conv-discrete}), almost surely, there exists $n_{2}$
such that, for $n\geq n_{2}$, $|\gamma_{T_{n+1}}(f)-\gamma_{T_{n}}(f)|<\iota$.
For $n\geq n_{1}\vee n_{2}$, we can then write:
\begin{multline}
\sum_{u\in\mathcal{U}_{\epsilon_{n}}^{(2)}}\xi_{u}\,\,\,\text{(by Equations (\ref{eq:terme-a}), (\ref{eq:terme-b}), (\ref{eq:terme-c}))}\\
\leq\frac{n+1}{n(1-e^{-a})}\left(\iota+\sum_{u\in\mathcal{U}_{\epsilon_{n}}^{(2)}}\sum_{v\in\boldsymbol{d}(u)}\xi_{v}\frac{1}{n+1}+\sum_{u\in\mathcal{U}_{\epsilon_{n}}^{(1)}}\xi_{u}\frac{1}{n}\right)\\
\text{(by Equation (\ref{eq:d-subset-a}))}\leq\frac{n+1}{n(1-e^{-a})}\left(\iota+\frac{1}{n}\right)\,.\label{eq:maj_U-2}
\end{multline}
Let $n\geq n_{1}\vee n_{2}$ and $t$ in $(T_{n},T_{n+1})$. We can
decompose
\begin{equation}
\mathcal{U}_{\epsilon_{n}}=\mathcal{U}_{\epsilon_{n}}^{(1)}(t)\sqcup\mathcal{U}_{\epsilon_{n}}^{(2)}(t)\text{ where }\mathcal{U}_{\epsilon_{n}}^{(1)}(t)=\mathcal{U}_{\epsilon_{n}}\cap\mathcal{U}_{e^{-t}}=\mathcal{U}_{e^{-t}}\,,\,\mathcal{U}_{\epsilon_{n}}^{(2)}(t)=\mathcal{U}_{\epsilon_{n}}\backslash\mathcal{U}_{e^{-t}}\,.\label{eq:dec-U-epsilon}
\end{equation}
For $u$ in $\mathcal{U}_{\epsilon_{n}}\backslash\mathcal{U}_{\epsilon_{n}}^{(1)}(t)$,
we set $\boldsymbol{d}(u,t)=\{v\in\mathcal{U}_{e^{-t}}\,:\,u=\mathbf{m}(v)\}$.
As $n\geq n_{1}$, $\mathbf{d}(u,t)=\mathbf{d}(u)$ and we have
\begin{equation}
\sum_{v\in\boldsymbol{d}(u,t)}\xi_{v}=\xi_{u}\,.\label{eq:d-subset-b}
\end{equation}
Similar to Equation (\ref{eq:dec-U-epsilon-02}), we have
\begin{equation}
\sum_{u\in\mathcal{U}_{e^{-t}}}\xi_{u}f(e^{t}\xi_{u})=\sum_{u\in\mathcal{U}_{\epsilon_{n}}^{(1)}(t)}\xi_{u}f(e^{t}\xi_{u})+\sum_{u\in\mathcal{U}_{\epsilon_{n}}^{(2)}(t)}\sum_{v\in\boldsymbol{d}(u,t)}\xi_{v}f(e^{t}\xi_{v})\,.\label{eq:dec-U-epsilon-03}
\end{equation}
 We fix $f$ continuous from $[0,1]$ to $\mathcal{\R}$, there exists
$n_{3}\in\N^{*}$ such that, for all $x,y\in[0,1]$, $|x-y|\leq1/n_{3}\Rightarrow|f(x)-f(y)|<\iota$.
Suppose that $n\geq n_{1}\vee n_{2}\vee n_{3}$. Then, using Equation
(\ref{eq:dec-U-epsilon}) and Equation (\ref{eq:d-subset-b}), we
have (for all $t\in[T_{n},T_{n+1}]$),
\begin{multline}
|\gamma_{t}(f)-\gamma_{T_{n}}(f)|=\left|\sum_{u\in\mathcal{U}_{e^{-t}}}\xi_{u}f(e^{t}\xi_{u})-\sum_{u\in\mathcal{U}_{\epsilon_{n}}^{(1)}(t)}\xi_{u}f(n\xi_{u})-\sum_{u\in\mathcal{U}_{\epsilon_{n}}^{(2)}(t)}\xi_{u}f(n\xi_{u})\right|\\
=\left|\sum_{u\in\mathcal{U}_{\epsilon_{n}}^{(1)}(t)}\xi_{u}f(e^{t}\xi_{u})+\sum_{u\in\mathcal{U}_{\epsilon_{n}}^{(2)}(t)}\sum_{v\in\boldsymbol{d}(u,t)}\xi_{v}f(e^{t}\xi_{v})-\sum_{u\in\mathcal{U}_{\epsilon_{n}}^{(1)}(t)}\xi_{u}f(n\xi_{u})-\sum_{u\in\mathcal{U}_{\epsilon_{n}}^{(2)}(t)}\xi_{u}f(n\xi_{u})\right|\\
\leq\left|\sum_{u\in\mathcal{U}_{\epsilon_{n}}^{(1)}(t)}\xi_{u}f(e^{t}\xi_{u})-\sum_{u\in\mathcal{U}_{\epsilon_{n}}^{(1)}(t)}\xi_{u}f(n\xi_{u})\right|+2\sum_{u\in\mathcal{U}_{\epsilon_{n}}^{(2)}(t)}\xi_{u}\Vert f\Vert_{\infty}\\
\leq\sum_{u\in\mathcal{U}_{\epsilon_{n}}^{(1)}(t)}\xi_{u}\iota+2\sum_{u\in\mathcal{U}_{\epsilon_{n}}^{(2)}(t)}\xi_{u}\Vert f\Vert_{\infty}\\
\text{(using Equation (\ref{eq:maj_U-2}) and as \ensuremath{\mathcal{U}_{\epsilon_{n}}^{(2)}(t)\subset\mathcal{U}_{\epsilon_{n}}^{(2)}})}\leq\iota+2\Vert f\Vert_{\infty}\frac{n+1}{n(1-e^{-a})}\left(\iota+\frac{1}{n}\right)\,.\label{eq:maj-temps-continu}
\end{multline}
Equations (\ref{eq:conv-discrete}) and (\ref{eq:maj-temps-continu})
prove the desired result.
\end{proof}
The set $\Bsym(1)$ is defined in Section \ref{subsec:Notations-1}.
\begin{thm}
[Central-limit Theorem] \label{thm:central-limit}We remember we
have Fact \ref{fact:theta} and we are under Hypothesis \ref{hyp:A},
\ref{hyp:conservative}, \ref{hyp:delta-step}, \ref{hyp:queue-pi}.
Let $q$ be in $\N^{*}$. For functions $f_{1}$, \ldots , $f_{q}$
which are continuous and in $\Bsym(1)$, we have
\[
\epsilon^{-q/2}(\gamma_{T}(f_{1}),\dots,\gamma_{T}(f_{q}))\underset{T\rightarrow+\infty}{\overset{\text{law}}{\longrightarrow}}\mathcal{N}(0,(K(f_{i},f_{j}))_{1\leq i,j\leq q})\,(\underline{\epsilon=e^{-T}})
\]
($K$ is given in Equation (\ref{eq:def-K})).
\end{thm}

\begin{proof}
Let $f_{1}$, \ldots , $f_{q}$ $\Bsym(1)$ and $v_{1},\dots,v_{q}\in\R$.
\\
\uline{First}, we develop the product below (remember that for
$u$ in $\mathcal{U}_{\epsilon}$, $\xi_{u}/\epsilon<1$ a.s.)
\begin{multline*}
\prod_{u\in\mathcal{U}_{\epsilon}}\left(1+\sqrt{\epsilon}\frac{\xi_{u}}{\epsilon}(iv_{1}f_{1}+\dots+iv_{q}f_{q})\left(\frac{\xi_{u}}{\epsilon}\right)\right)=\\
\exp\left(\sum_{u\in\mathcal{U}_{\epsilon}}\log\left[1+\sqrt{\epsilon}\Id\times(iv_{1}f_{1}+\dots+iv_{q}f_{q})\left(\frac{\xi_{u}}{\epsilon}\right)\right]\right)=\\
\text{(for }\epsilon\text{ small enough)}\\
\exp\left(\sum_{u\in\mathcal{U_{\epsilon}}}\sum_{k\geq1}\frac{(-1)^{k+1}}{k}\epsilon^{k/2}(\Id\times(iv_{1}f_{1}+\dots+iv_{q}f_{q}))^{k}\left(\frac{\xi_{u}}{\epsilon}\right)\right)=\\
\exp\left(\frac{1}{\sqrt{\epsilon}}\gamma_{T}(iv_{1}f_{1}+\dots+iv_{q}f_{q})+\frac{1}{2}\gamma_{T}(\Id\times(v_{1}f_{1}+\dots+v_{q}f_{q})^{2})+R_{\epsilon}\right)\,,
\end{multline*}
where 
\begin{eqnarray*}
R_{\epsilon} & = & \sum_{k\geq3}\sum_{u\in\mathcal{U}_{\epsilon}}\frac{(-1)^{k+1}}{k}\epsilon^{k/2-1}\xi_{u}\left(\frac{\xi_{u}}{\epsilon}\right)^{k-1}(iv_{1}f_{1}+\dots+iv_{q}f_{q})^{k}\left(\frac{\xi_{u}}{\epsilon}\right)\\
 & = & \sum_{k\geq3}\frac{(-1)^{k+1}}{k}\epsilon^{k/2-1}\gamma_{T}((\Id)^{k-1}(iv_{1}f_{1}+\dots+iv_{q}f_{q})^{k})\,,
\end{eqnarray*}
\[
|R_{\epsilon}|\leq\sum_{k\geq3}\frac{\epsilon^{k/2-1}}{k}(|v_{1}|\Vert f_{1}\Vert_{\infty}+\dots+|v_{q}|\Vert f_{q}\Vert_{\infty})^{k}=O(\sqrt{\epsilon})\,.
\]
We have, for some constant $C$, (we use that: $x\in\R\Rightarrow|e^{ix}|=1$)
\begin{multline*}
\E\left(\left|\exp\left(\frac{1}{\sqrt{\epsilon}}\gamma_{T}(iv_{1}f_{1}+\dots+iv_{q}f_{q})+\frac{1}{2}\gamma_{T}(\Id\times(v_{1}f_{1}+\dots+v_{q}f_{q})^{2})+R_{\epsilon}\right)\right.\right.\\
\left.\left.-\exp\left(\frac{1}{\sqrt{\epsilon}}\gamma_{T}(iv_{1}f_{1}+\dots+iv_{q}f_{q})+\frac{1}{2}\eta(\Phi(\Id\times(v_{1}f_{1}+\dots+v_{q}f_{q})^{2})\right)\right|\right)\\
\leq\E\left(C\left|\frac{1}{2}\gamma_{T}(\Id\times(v_{1}f_{1}+\dots+v_{q}f_{q})^{2})-\frac{1}{2}\eta(\Phi(\Id\times(v_{1}f_{1}+\dots+v_{q}f_{q})^{2})+R_{\epsilon}\right|\right)\\
\text{(by Proposition \ref{prop:convergence-ps})}\underset{\epsilon\rightarrow0}{\longrightarrow}0\,.
\end{multline*}

\uline{Second}, we develop the same product in a different manner.
We have (the order on $\mathcal{U}$ is defined in Section \ref{subsec:Fragmentation-chains})
\begin{multline*}
\prod_{u\in\mathcal{U}_{\epsilon}}\left(1+\sqrt{\epsilon}\frac{\xi_{u}}{\epsilon}(iv_{1}f_{1}+\dots+iv_{q}f_{q})\left(\frac{\xi_{u}}{\epsilon}\right)\right)=\\
\sum_{k\geq0}\epsilon^{-k/2}i^{k}\sum_{1\leq j_{1},\dots,j_{k}\leq q}v_{j_{1}}\dots v_{j_{k}}\sum_{\begin{array}{c}
u_{1},\dots,u_{k}\in\mathcal{U_{\epsilon}}\\
u_{1}<\dots<u_{k}
\end{array}}\xi_{u_{1}\dots}\xi_{u_{k}}f_{j_{1}}\left(\frac{\xi_{u_{1}}}{\epsilon}\right)\dots f_{j_{k}}\left(\frac{\xi_{u_{k}}}{\epsilon}\right)=\\
\text{(a detailed proof can be found in Section \ref{subsec:Detailed-proof-of-2})}\\
\sum_{k\geq0}\epsilon^{-k/2}i^{k}\sum_{1\leq j_{1},\dots,j_{k}\leq q}v_{j_{1}}\dots v_{j_{k}}\frac{1}{k!}\gamma_{T}^{\odot k}(f_{j_{1}}\otimes\dots\otimes f_{j_{k}})\,.
\end{multline*}
We have, for all $k$,
\begin{multline*}
\left|\epsilon^{-k/2}\sum_{1\leq j_{1},\dots,j_{k}\leq q}v_{j_{1}}\dots v_{j_{k}}\frac{1}{k!}\E(\gamma_{T}^{\odot k}(f_{j_{1}}\otimes\dots\otimes f_{j_{k}}))\right|\\
\leq\epsilon^{-k/2}\times\frac{q^{k}\sup(|v_{1}|,\dots,|v_{q}|)^{k}\sup(\Vert f_{1}\Vert_{\infty},\dots,\Vert f_{q}\Vert_{\infty})^{k}}{k!}.
\end{multline*}
 So, by Corollary \ref{cor:lim-case-q-odd}, Proposition \ref{prop:conv-q-pair}
and Equation (\ref{eq:lien-odot-B}), we get that 
\begin{multline*}
\E\left(\prod_{u\in\mathcal{U}_{\epsilon}}\left(1+\sqrt{\epsilon}\frac{\xi_{u}}{\epsilon}(iv_{1}f_{1}+\dots+iv_{q}f_{q})\left(\frac{\xi_{u}}{\epsilon}\right)\right)\right)\\
\underset{\epsilon\rightarrow0}{\longrightarrow}\sum_{\begin{array}{c}
k\geq0\\
k\text{ even}
\end{array}}(-1)^{k/2}\sum_{1\leq j_{1},\dots,j_{k}\leq q}\frac{1}{k!}\sum_{I\in I_{k}}\prod_{\{a,b\}\in I}V(v_{j_{a}}f_{j_{a}},v_{j_{b}}f_{j_{b}})\\
\text{(a detailed proof can be found in Section \ref{subsec:Detailed-proof-of-1})}\\
\text{(using Equation (\ref{eq:card-I-q}))}=\sum_{\begin{array}{c}
k\geq0\\
k\text{ even}
\end{array}}\frac{(-1)^{k/2}}{2^{k/2}(k/2)!}\sum_{1\leq j_{1},\dots,j_{k}\leq q}V(v_{j_{1}}f_{j_{1}},v_{j_{2}}f_{j_{2}})\dots V(v_{j_{k-1}}f_{j_{k-1}},v_{j_{k}}f_{j_{k}})\\
=\sum_{\begin{array}{c}
k\geq0\\
k\text{ even}
\end{array}}\frac{(-1)^{k/2}}{2^{k/2}(k/2)!}\left(\sum_{1\leq j_{1},j_{2}\leq q}v_{j_{1}}v_{j_{2}}V(f_{j_{1}},f_{j_{2}})\right)^{k/2}\\
=\exp\left(-\frac{1}{2}\sum_{1\leq j_{1},j_{2}\leq q}v_{j_{1}}v_{j_{2}}V(f_{j_{1}},f_{j_{2}})\right)\,.
\end{multline*}
\uline{In conclusion}, we have
\begin{multline*}
\E\left(\exp\left(\frac{1}{\sqrt{\epsilon}}\gamma_{T}(iv_{1}f_{1}+\dots+iv_{q}f_{q})\right)\right)\\
\underset{\epsilon\rightarrow0}{\longrightarrow}\exp\left(-\frac{1}{2}\eta(\Phi(\Id\times(v_{1}f_{1}+\dots+v_{q}f_{q})^{2}))-\frac{1}{2}\sum_{1\leq j_{1},j_{2}\leq q}v_{j_{1}}v_{j_{2}}V(f_{j_{1}},f_{j_{2}})\right)\,.
\end{multline*}
So we get the desired result with, for all $f$, $g$,
\begin{equation}
K(f,g)=\eta(\Phi(\Id\times fg)+V(f,g))\label{eq:def-K}
\end{equation}
($V$ is defined in Equation (\ref{eq:def-V})).
\end{proof}
\bibliographystyle{amsalpha}
\bibliography{biblio-fragmentation}

\section{Appendix\label{sec:Appendix}}

\subsection{Detailed proof\label{subsec:Detailed-proof-of} of a bound appearing
in the proof of Lemma \ref{lem:conv-nouille-libre}}
\begin{lem}
We have (for any $f$ appearing in the proof of Lemma \ref{lem:conv-nouille-libre})
\[
\E(\1_{A_{u}=f(u),\forall u\in\mathcal{T}_{2}}|\mathcal{L}_{2},\mathcal{T}_{2},m_{2})\leq\prod_{u\in\mathcal{T}_{2}\backslash\{0\}}e^{-(\#f(u)-1)(T_{u}-T_{\boldsymbol{m}(u)})}\,.
\]
\end{lem}

\begin{proof}
We want to show this by recurrence on the cardinality of $\mathcal{T}_{2}$.

If $\#\mathcal{T}_{2}=1$, then $\mathcal{T}_{2}=\{0\}$ and the claim
is true.

Suppose now that $\#\mathcal{T}_{2}=k$ and the claim is true up to
the cardinality $k-1$. There exists $v$ in $\mathcal{T}_{2}$ such
that $(v,i)$ is not in $\mathcal{T}_{2}$, for any $i$ in $\N^{*}$.
We set $\mathcal{T}_{2}'=\mathcal{T}_{2}\backslash\{v\}$, $\mathcal{L}_{2}'=\mathcal{L}_{2}\backslash\{v\}$,
$m_{2}'\,:\,u\in\mathcal{T}_{2}'\rightarrow(\xi_{u},\inf\{i,i\in A_{u}\})$.
We set $f(v)=\{i_{1},\dots,i_{p}\}$ (with $i_{1}<\dots<i_{p}$),
$f(\boldsymbol{m}(v))=\{i_{1},\dots,i_{p},i_{p+1},\dots,i_{q}\}$
(with $i_{p+1}<\dots<i_{q}$). We suppose $m_{2}(v)=(\xi_{v},i_{1})$
because if $m_{2}(v)=(\xi_{v},j)$ with $j\neq i_{1}$ then $A_{v}\neq f(v)$
for all $\omega$, and then the left-hand side of the inequality above
is zero. We have 
\begin{multline*}
\E(\1_{A_{u}=f(u),\forall u\in\mathcal{T}_{2}}|\mathcal{L}_{2},\mathcal{T}_{2},m_{2})\\
=\E(\1_{A_{u}=f(u),\forall u\in\mathcal{T}_{2}'}\E(\1_{A_{v=f(v)}}|\mathcal{L}_{2},\mathcal{T}_{2},m_{2},(A_{u},u\in\mathcal{T}_{2}'))|\mathcal{L}_{2},\mathcal{T}_{2},m_{2})\\
=\E(\1_{A_{u}=f(u),\forall u\in\mathcal{T}_{2}'}\E(\1_{i_{1},\dots,i_{p}\in A_{v}}\1_{i_{p+1},\dots,i_{q}\notin A_{v}}|\mathcal{L}_{2},\mathcal{T}_{2},m_{2},(A_{u},u\in\mathcal{T}_{2}'))|\mathcal{L}_{2},\mathcal{T}_{2},m_{2})\\
\text{(remember we condition on \ensuremath{m_{2}}, so the \ensuremath{\1_{i_{1},\dots,i_{p}}}can be replaced by \ensuremath{\1_{i_{2},\dots,i_{p}}})}\\
=\E(\1_{A_{u}=f(u),\forall u\in\mathcal{T}_{2}'}\E(\1_{i_{2},\dots,i_{p}\in A_{v}}\1_{i_{p+1},\dots,i_{q}\notin A_{v}}|\mathcal{L}_{2},\mathcal{T}_{2},m_{2},(A_{u},u\in\mathcal{T}_{2}'))|\mathcal{L}_{2},\mathcal{T}_{2},m_{2})\\
\leq\E(\1_{A_{u}=f(u),\forall u\in\mathcal{T}_{2}'}\E(\1_{i_{2},\dots,i_{p}\in A_{v}}|\mathcal{L}_{2},\mathcal{T}_{2},m_{2},(A_{u},u\in\mathcal{T}_{2}'))|\mathcal{L}_{2},\mathcal{T}_{2},m_{2})\\
\text{(because the \ensuremath{(Y_{j})} introduced in Section \ref{subsec:Tagged-fragments} are independant) }\\
=\E(\1_{A_{u}=f(u),\forall u\in\mathcal{T}_{2}'}\prod_{r=2}^{p}\E(\1_{i_{r}\in A_{v}}|\mathcal{L}_{2},\mathcal{T}_{2},m_{2},(A_{u},u\in\mathcal{T}_{2}'))|\mathcal{L}_{2},\mathcal{T}_{2},m_{2})\\
\text{(because of Equation (\ref{eq:proba-rester-dans-fragment}))}\\
\text{(if \ensuremath{v\in\mathcal{L}_{2}} then \ensuremath{\prod_{r=2}^{p}\dots} is empty and thus \ensuremath{=1})}\\
=\E(\1_{A_{u}=f(u),\forall u\in\mathcal{T}_{2}'}\prod_{r=2}^{p}\widetilde{\xi}_{v}|\mathcal{L}_{2},\mathcal{T}_{2},m_{2})\\
\text{(by Equation (\ref{eq:def-T_u}) and Proposition \ref{prop:(reformulation-of-Proposition})}\\
=e^{-(\#f(v)-1)(T_{v}-T_{\boldsymbol{m}(v)})}\E(\1_{A_{u}=f(u),\forall u\in\mathcal{T}_{2}'}|\mathcal{L}_{2},\mathcal{T}_{2},m_{2})\\
=e^{-(\#f(v)-1)(T_{v}-T_{\boldsymbol{m}(v)})}\E(\E(\1_{A_{u}=f(u),\forall u\in\mathcal{T}_{2}'}|\mathcal{L}_{2}',\mathcal{T}_{2}',m_{2}')|\mathcal{L}_{2},\mathcal{T}_{2},m_{2})\\
\text{(by recurrence)}\\
\leq\prod_{u\in\mathcal{T}_{2}\backslash\{0\}}e^{-(\#f(u)-1)(T_{u}-T_{\boldsymbol{m}(u)})}\,.
\end{multline*}
 
\end{proof}

\subsection{\label{subsec:Detailed-proof-of-3}Detailed proof of a bound appearing
in the proof of Corollary \ref{cor:maj-fonction-centree}}
\begin{lem}
Let $q$ be in $\N$, we have
\[
\sum_{k'\in[q]}1+(q-k'-1)_{+}\leq(q+1)^{2}\,.
\]
\end{lem}

\begin{proof}
We have 
\begin{eqnarray*}
\sum_{k'\in[q]}1+(q-k'-1)_{+} & = & q+\sum_{k'\in[q-2]}(q-k'-1)\\
 & = & q+\sum_{i=1}^{q-2}i\\
 & \leq & \frac{q(q+1)}{2}\\
 & \leq & (q+1)^{2}\,.
\end{eqnarray*}
\end{proof}

\subsection{\label{subsec:Detailed-proof-of-1}Detailed proof of an equality
appearing in the proof of Theorem \ref{thm:central-limit}}
\begin{lem}
Let $q\in\N^{*}$. Suppose we have $q$ functions $g_{1}$, \dots ,
$g_{q}$ in $\Bsym(1)$. Then, for all $k$ even ($k$ in $\N$)
\begin{equation}
\sum_{1\leq j_{1},\dots,j_{k}\leq q}\sum_{I\in\mathcal{I}_{k}}\prod_{\{a,b\}\in I}V(g_{j_{a}},g_{j_{b}})=\frac{k!}{2^{k/2}(k/2)!}\sum_{1\leq j_{1},\dots,j_{k}\leq q}V(g_{j_{1}},g_{j_{2}})\dots V(g_{j_{k-1}},g_{j_{k}})\,.\label{eq:denombrement-01}
\end{equation}
\end{lem}

\begin{proof}
We set 
\[
\sum_{1\leq j_{1},\dots,j_{k}\leq q}\sum_{I\in\mathcal{I}_{k}}\prod_{\{a,b\}\in I}V(g_{j_{a}},g_{j_{b}})=(1)\,,
\]
\[
\sum_{1\leq j_{1},\dots,j_{k}\leq q}V(g_{j_{1}},g_{j_{2}})\dots V(g_{j_{k-1}},g_{j_{k}})=(2)\,.
\]
Suppose, for some $k$, we have $i_{1}$, \dots , $i_{k}$ in $[q]$,
distinct. There exists $N_{1}$, $N_{2}$ such that:
\begin{itemize}
\item the term (1) has $N_{1}$ terms $V(g_{i_{1}},g_{i_{2}})\dots V(g_{i_{k-1}},g_{i_{k}})$
(up to permutations, that is we consider that $V(g_{i_{3}},g_{i_{4}})V(g_{i_{2}},g_{i_{1}})\dots V(g_{i_{k-1}},g_{i_{k}})$
and $V(g_{i_{1}},g_{i_{2}})\dots V(g_{i_{k-1}},g_{i_{k}})$ are the
same term),
\item the term (2) has $N_{2}$ terms $V(g_{i_{1}},g_{i_{2}})\dots V(g_{i_{k-1}},g_{i_{k}})$
(again, up to permutations).
\end{itemize}
These numbers $N_{1}$, $N_{2}$ do not depend on $i_{1}$, \dots ,
$i_{k}$. In the case where the indexes $i_{1}$, \dots , $i_{k}$
are not distinct, we can find easily the number of terms equal to
$V(g_{i_{1}},g_{i_{2}})\dots V(g_{i_{k-1}},g_{i_{k}})$ in terms (1),
(2). For example, if $i_{2}=i_{1}$ and $i_{1},\,i_{3},\dots,i_{k}$
are distinct, then 
\begin{itemize}
\item the term (1) has $2N_{1}$ terms $V(g_{i_{1}},g_{i_{2}})\dots V(g_{i_{k-1}},g_{i_{k}})$,
\item the term (2) has $2N_{2}$ terms $V(g_{i_{1}},g_{i_{2}})\dots V(g_{i_{k-1}},g_{i_{k}})$ 
\end{itemize}
(we multiply simply by the number of $\sigma$ in $\mathcal{S}_{k}$
such that $(i_{1},i_{2},\dots,i_{k})=(i_{\sigma(1)},i_{\sigma(2)},\dots,i_{\sigma(k)})$).
We do not need to know $N_{1}$ and $N_{2}$ but we need to know $N_{1}/N_{2}$.
By taking $V(g,f)$ to be $1$ for all $g$, $f$, we see that $N_{1}/N_{2}=\#\mathcal{I}_{k}=k!/(2^{k/2}(k/2)!)$. 
\end{proof}

\subsection{\label{subsec:Detailed-proof-of-2}Detailed proof of an equality
appearing in the proof of Theorem \ref{thm:central-limit}}
\begin{lem}
Let $f_{1}$, \ldots , $f_{q}$ $\Bsym(1)$ , $k$ in $\N$ and $v_{1},\dots,v_{q}\in\R$.
We have
\begin{multline*}
\sum_{1\leq j_{1},\dots,j_{k}\leq q}v_{j_{1}}\dots v_{j_{k}}\gamma_{T}^{\odot k}(f_{j_{1}}\otimes\dots\otimes f_{j_{k}})\\
=k!\sum_{1\leq j_{1},\dots,j_{k}\leq q}v_{j_{1}}\dots v_{j_{k}}\sum_{\begin{array}{c}
u_{1},\dots,u_{k}\in\mathcal{U_{\epsilon}}\\
u_{1}<\dots<u_{k}
\end{array}}\xi_{u_{1}\dots}\xi_{u_{k}}f_{j_{1}}\left(\frac{\xi_{u_{1}}}{\epsilon}\right)\dots f_{j_{k}}\left(\frac{\xi_{u_{k}}}{\epsilon}\right)\,.
\end{multline*}
\end{lem}

\begin{proof}
We have 
\begin{multline*}
\sum_{1\leq j_{1},\dots,j_{k}\leq q}v_{j_{1}}\dots v_{j_{k}}\gamma_{T}^{\odot k}(f_{j_{1}}\otimes\dots\otimes f_{j_{k}})\\
=\sum_{1\leq j_{1},\dots,j_{k}\leq q}v_{j_{1}}\dots v_{j_{k}}\sum_{a:[k]\hookrightarrow\mathcal{U}_{\epsilon}}\xi_{a(1)}\dots\xi_{a(k)}f_{j_{1}}\left(\frac{\xi_{a(1)}}{\epsilon}\right)\dots f_{j_{k}}\left(\frac{\xi_{a(k)}}{\epsilon}\right)\\
\text{(for all injection }a\text{, there is exactly one }\sigma_{a}\in\mathcal{S}_{k}\text{ such that }a(\sigma_{a}(1))<\dots<a(\sigma_{a}(k))\text{)}\\
=\sum_{1\leq j_{1},\dots,j_{k}\leq q}v_{j_{1}}\dots v_{j_{k}}\sum_{a:[k]\hookrightarrow\mathcal{U}_{\epsilon}}\xi_{a(\sigma_{a}(1))}\dots\xi_{a(\sigma_{a}(k))}f_{j_{\sigma_{a}(1)}}\left(\frac{\xi_{a(\sigma_{a}(1))}}{\epsilon}\right)\dots f_{j_{\sigma_{a}(k)}}\left(\frac{\xi_{a(\sigma_{a}(k))}}{\epsilon}\right)\\
\text{(for }\tau\in\mathcal{S}_{k}\text{, we set }\mathcal{E}(\tau)=\{a:[k]\hookrightarrow\mathcal{U}_{\epsilon}\,:\,\sigma_{a}=\tau\}\text{)}\\
=\sum_{1\leq j_{1},\dots,j_{k}\leq q}v_{j_{1}}\dots v_{j_{k}}\sum_{\tau\in\mathcal{S}_{k}}\sum_{a\in\mathcal{E}(\tau)}\xi_{a(\tau(1))}\dots\xi_{a(\tau(k))}f_{j_{\tau(1)}}\left(\frac{\xi_{a(\tau(1))}}{\epsilon}\right)\dots f_{j_{\tau(k)}}\left(\frac{\xi_{a(\tau(k))}}{\epsilon}\right)\\
=\sum_{\tau\in\mathcal{S}_{k}}\sum_{1\leq j_{1},\dots,j_{k}\leq q}v_{j_{1}}\dots v_{j_{k}}\sum_{a\in\mathcal{E}(\tau)}\xi_{a(\tau(1))}\dots\xi_{a(\tau(k))}f_{j_{\tau(1)}}\left(\frac{\xi_{a(\tau(1))}}{\epsilon}\right)\dots f_{j_{\tau(k)}}\left(\frac{\xi_{a(\tau(k))}}{\epsilon}\right)\\
=\sum_{\tau\in\mathcal{S}_{k}}\sum_{1\leq j_{1},\dots,j_{k}\leq q}v_{j_{\tau(1)}}\dots v_{j_{\tau(k)}}\sum_{a\in\mathcal{E}(\tau)}\xi_{a(\tau(1))}\dots\xi_{a(\tau(k))}f_{j_{\tau(1)}}\left(\frac{\xi_{a(\tau(1))}}{\epsilon}\right)\dots f_{j_{\tau(k)}}\left(\frac{\xi_{a(\tau(k))}}{\epsilon}\right)\\
=\sum_{\tau\in\mathcal{S}_{k}}\sum_{1\leq j_{1},\dots,j_{k}\leq q}v_{j_{\tau(1)}}\dots v_{j_{\tau(k)}}\sum_{\overset{u_{1},\dots,u_{k}\in\mathcal{U}_{\epsilon}}{u_{1}<\dots<u_{k}}}\xi_{u_{1}}\dots\xi_{u_{k}}f_{j_{\tau(1)}}\left(\frac{\xi_{u_{1}}}{\epsilon}\right)\dots f_{j_{\tau(k)}}\left(\frac{\xi_{u_{k}}}{\epsilon}\right)
\end{multline*}
The application (``$\hookrightarrow$'' means that an application
is injective)
\[
(a\,:\,[k]\rightarrow[q]\,,\,\tau\,:\,[k]\hookrightarrow[k])\overset{\Theta}{\longrightarrow}a\circ\tau\,
\]
is such that
\[
\forall b\,:\,[k]\rightarrow[q]\,,\,\#\Theta^{-1}(\{b\})=k!\,.
\]
So the above quantity is equal to 
\[
k!\sum_{1\leq j_{1},\dots,j_{k}\leq q}v_{j_{1}}\dots v_{j_{k}}\sum_{\overset{u_{1},\dots,u_{k}\in\mathcal{U}_{\epsilon}}{u_{1}<\dots<u_{k}}}\xi_{u_{1}}\dots\xi_{u_{k}}f_{j_{1}}\left(\frac{\xi_{u_{1}}}{\epsilon}\right)\dots f_{j_{k}}\left(\frac{\xi_{u_{k}}}{\epsilon}\right)\,.
\]
\end{proof}

\end{document}